\documentclass[12pt]{article}

\usepackage[left=1in,right=1in,top=1in,bottom=1in]{geometry}

\usepackage{subfigure}
\usepackage{graphicx}
\usepackage{amsfonts}
\usepackage{url}
\usepackage{amsmath,amssymb}
\usepackage{enumerate}
\usepackage{amsthm,amsmath}
\usepackage{float}

\usepackage{mathtools}

\usepackage[colorlinks]{hyperref}

\newtheorem{algorithm}{Algorithm}
\newtheorem{remark}{Remark}

\theoremstyle{definition}
\newtheorem{model}{Example Model}
\newtheorem{example}{Example}

\newcommand{\R}{{\mathbb R}}
\newcommand{\Z}{{\mathbb Z}}

\newcommand{\Ql}{Q}

\newcommand{\hl}{h_{\ell}}

\newcommand\qedex{\hfill $\square$}

\title{Low variance couplings for stochastic models of intracellular processes with time-dependent rate functions}

\author{
David F. Anderson\thanks{Department of Mathematics, University of
  Wisconsin, Madison, USA.  anderson@math.wisc.edu, grant support from NSF-DMS-1318832 and Army Research Office grant W911NF-14-1-0401.},
\and
Chaojie Yuan\thanks{Department of Mathematics, University of
  Wisconsin, Madison, USA.  cyuan25@math.wisc.edu.
}}

\begin{document}

\maketitle

\begin{abstract}
 	A number of coupling strategies are presented for stochastically modeled biochemical processes with time-dependent parameters.  In particular, the \textit{stacked} coupling is introduced and is shown via a number of examples to provide an exceptionally low variance between the generated paths.  This coupling will be useful in the numerical computation of parametric sensitivities and the fast estimation of expectations via multilevel Monte Carlo methods.  We provide the requisite estimators in both cases.
\end{abstract}

\section{Introduction}
\label{sec:intro}

We consider stochastic models of intracellular processes whose rate functions, perhaps due to changes in temperature, volume, voltage (in the case of neural networks), or even some external (possibly random) forcing, depend explicitly on time.  Specifically, we consider the standard discrete-space, continuous-time Markov chain model, typically simulated via Gillespie's algorithm \cite{gillespie1976general,gillespie1977exact} or the next reaction method \cite{anderson2007modified,Gibson2000}, whose propensity functions take the form $\lambda_k(t,X(t))$, with $X(t)$ being the vector whose $i$th component gives the count of species $i$ at time $t$, and the parameter $k$ enumerates over the reaction channels.

We provide an example model to solidify notation.

 \begin{model}\label{example:first}
Consider the following standard model of  transcription and translation,
\begin{align*}
	\emptyset &\overset{60}\to  M \qquad \qquad	M \overset{100}\to  M + P\qquad\qquad M \overset{1}{\to} \emptyset \qquad\qquad P \overset{1}{\to} \emptyset,
\end{align*}
where $M$ represents mRNA and $P$ represents proteins.  The notation $M\overset{100}\to M+P$ implies:
\begin{enumerate}[(i)]
\item at least one mRNA molecule is required to be present for the reaction to take place (implied by the $M$ on the left of the arrow), and the net change when the reaction takes place is the addition of one protein molecule, and 

\item the mass-action rate parameter is 100, yielding a propensity function of $\lambda_2(x) = 100 x_1$ (where we arbitrarily call $M$ the first species and $P$, the protein, the second).
\end{enumerate}
This model (though not necessarily the specific choice of rate constants) is one of the simplest and most popular models in the study of cell biology.
 Note that the rate parameters (60, 100, 1, and 1) are fixed constants in the model.  This need not be the case.  For example, perhaps the system is subject to dark-light cycles that oscillate over a 24 hour time-period.  In this case, a better model may be to change the propensity function of the first reaction, which is currently $\lambda_1(x) = 60$, to 
\[
	\lambda_1(t,x)=60+15\sin\left(\tfrac{2\pi t}{24}\right).
\]
Of course, it could be that the other rates change explicitly in time as well. \qedex
\end{model}

 In this situation of time dependent propensity functions, the basic versions of the Gillespie algorithm and the next reaction method can  be quite slow in generating sample paths for even simple models.    The slowdown occurs because the calculation for the time required before  the occurrence of the next event in the system becomes a hitting time problem of the following form:  find $\Delta$ solving
 \begin{equation}\label{eq:int_to_estimate1}
 	\int_t^{t+\Delta} g(s) ds = \mathcal{E},
 \end{equation}
 where $g$ is some function (depending on the time-varying parameters), and $\mathcal{E}$ is some positive value (often a unit exponential random variable).  Solving for $\Delta$ typically requires a numerical estimation of the integral, thereby dramatically slowing the simulation, often by factors in the hundreds or thousands \cite{anderson2007modified,AET2015, lemaire2016exact,thanh2015simulation,voliotis2016stochastic}.
  A way around the necessity for estimating the integrals in \eqref{eq:int_to_estimate1} is to utilize thinning procedures \cite{GlynnAsmussen2007}. These  methods have been shown to dramatically reduce the cost of path simulation in the present context \cite{voliotis2016stochastic}. A similar idea using the rejection based stochastic simulation algorithm \cite{marchetti2016hrssa, thanh2015simulation} can be applied to save computation cost by delaying the update of propensity functions. Other applications of thinning have also appeared in the literature recently. Hybrid simulation schemes have been  proposed to simulate jump-diffusion processes that combine the jump process and diffusion approximation \cite{duncan2015hybrid}. Exact trajectories of a class of Piecewise deterministic Markov processes can also be simulated using thinning \cite{lemaire2016exact}.

 A particularly active, and fruitful, area of research over the last decade for models with time independent parameters has focused  on the following question:  given two mathematical models for a   system, which we will denote by $X$ and $Z$, respectively, and a functional $f$ of the model output, how can we generate the pair $(X,Z)$ so that $\text{Var}(f(X)-f(Z))$ is small? That is, how can we \textit{couple} the processes. The two most common applications of this research are (i) parametric sensitivity analysis, in which $X$ and $Z$ differ by the perturbation of some parameters of interest \cite{AndCFD2012, GuptaKhammash, Khammash2010, SAR2013}, and (ii) multilevel Monte Carlo where $X$ is an exact process and $Z$ is a tau-leap discretization, or both $X$ and $Z$ are tau-leap discretizations, though with different step-sizes \cite{AndHigham2012,AHS2014,AHS2016, AHS2017,MTV2016a}.  
   
In this paper, we  study the coupling question in the context of time-dependent parameters.   In particular, we provide  a variety of equivalent stochastic equations for the models under consideration and use these to develop a number of reasonable coupling strategies.  We discover that one coupling, here termed the \textit{stacked coupling}, is far  more efficient than the others on the examples we have considered. Specifically,  it produces pairs of paths with exceptionally low variance, often hundreds   of times lower than the variance between the paths produced by the other couplings.  We include the detailed construction of the other couplings in case they find use in  contexts not considered here.

The outline of the remainder of the paper is as follows.  In Section \ref{sec:math_models}, we will formally introduce the models considered in this paper.  While the stochastic equations are different (they are built on different probability spaces) each produces paths whose distributions agree with the desired model (similar to how the Gillespie algorithm and next reaction method are different, but produce equivalent paths  in distribution).    In Section \ref{sec:couplings1}, we provide multiple different couplings, which arise naturally via the different stochastic equations provided in Section \ref{sec:math_models}.  In Section \ref{sec:applications} we present applications and examples, focusing on parametric sensitivity analysis and multi-level Monte Carlo.    We end with some conclusions and future directions for research in Section \ref{sec:discussion}.

\section{Mathematical models and their (multiple) stochastic equation representations} 
 \label{sec:math_models}
 
%
 
 \subsection{The time homogeneous case}
 The most common mathematical model utilized for  intracellular biochemical processes in which the counts, or abundances, of the constituent species are tracked is a  time-homogeneous, discrete-space, continuous-time Markov process. We review this model first before transitioning to the non-homogeneous case.
 
 We suppose that the model includes $d$ chemical species whose abundances can change due to any of $K$ different possible reactions, where $K$ is a positive integer.    Let $X(t)$ be the vector in $\Z^d_{\ge0}$ whose $i$th component, $X_i(t)$, gives the abundance of  the $i$th species at time $t\ge 0$.  Let $\zeta_k \in \Z^d$ be the \textit{reaction vector} for the $k$th reaction channel, so that if reaction $k$ occurs at time $t^*$ we have
\[
	X(t^*) = X(t^*-) + \zeta_k,
\]
where $X(t^*-)$ denotes the state of the system just prior to time $t^*$.    Simple bookkeeping then says that if $R_k(t)$ is the number of times reaction $k \in \{1,\dots,K\}$ has occurred by time $t$, then
\begin{equation}\label{eq:counting_8976786}
	X(t) = X(0) + \sum_{k=1}^K R_k(t) \zeta_k.
\end{equation}
The $R_k$ are counting processes  and we assume that they have state dependent propensity functions (termed \textit{intensity functions} in the probability literature) $\lambda_k(\cdot)$ so that in the limit as $h\to 0$
\begin{align}
\begin{split}
	P(\text{reaction $k$ occurs in $[t,t+h)$}\ | \ \mathcal{F}_t) &= \lambda_k(X(t)) h + o(h)\\
	P(\text{two or more reactions occur in $[t,t+h)$}\ | \ \mathcal{F}_t) &= o(h)\\
	P(\text{no reaction occurs in $[t,t+h)$}\ | \ \mathcal{F}_t) &=1 - \sum_{k}\lambda_k(X(t)) h + o(h),
	\end{split}
	\label{eq:assump_jump_rates_homogeneous}
\end{align}
where $\mathcal{F}_t$ represents all the information learned about the process up to time $t$ ($\mathcal F_t$ is termed a \textit{filtration} in the probability literature), and $o(h)$ represents a term that is substantially smaller than $h$ (to be precise: it is a function, say $g$, for which  $g(h)/h \to 0$, as $h \to 0$).  

The process satisfying the assumptions above is a continuous time Markov chain in $\Z^d_{\ge 0}$, and there are multiple ways to specify it \cite{AK2011, AK2015}.  For example, it is the Markov chain on $\Z^d_{\ge 0}$ with master equation (termed the \textit{forward equation} in the probability literature)
\begin{equation*}
	\frac{d}{dt} P_\mu(x,t) = \sum_{k = 1}^K P_\mu(x-\zeta_k,t) \lambda_k(x-\zeta_k) - \sum_{k=1}^K \lambda_k(x)P_\mu(x,t),
\end{equation*}
where $P_\mu(x,t)$ is the probability $X(t) = x$ assuming an initial distribution of $\mu$.  The process may also be specified as the Markov chain with infinitesimal generator
\[
	Af(x) = \sum_k \lambda_k(x) (f(x+\zeta_k)-f(x)),
\]
for functions $f:\Z^d_{\ge 0} \to \R$.

Conversely, we can specify the model by specifying the counting processes in \eqref{eq:counting_8976786}. 
There are multiple ways to do this, depending on the probability space you choose to work with.  The different representations imply different algorithms, a critical relation at the heart of this paper.  Below is a non-exhaustive list of different possible representations.

\begin{enumerate}
\item The most widely used representation is usually called the  \textit{random time change representation}, and was developed by Thomas Kurtz \cite{Kurtz80}.  It has been utilized widely for both the development of computational methods \cite{anderson2007modified,Anderson2007b,AndCFD2012,AET2015,AndHigham2012,AndSkubak,GuptaKhammash,GK2014b,GK2014,MTV2016a}  and for analytical purposes \cite{ACK2016,AEJ2014,AndersonGangulyKurtz,AK2015}.  For this representation, we start with independent unit-rate Poisson processes $Y_k$ (one for reach reaction channel) and define the process $X$ as the solution to 
\begin{equation}\label{eq:RTC_main}
X(t) = X(0)  + \sum_{k=1}^K Y_k \big( \int_0^t \lambda_k(X(s)) ds\big) \zeta_k.
\end{equation}
Hence, in this case $R_k(t) = Y_k(\int_0^t \lambda_k(X(s)) ds )$. Exact simulation of this representation is equivalent to the next reaction method of Gibson and Bruck \cite{anderson2007modified,Gibson2000}.

\item The second representation begins with a space-time unit-rate Poisson point process $N$ (so the number of points, $N(A)$, in a region $A$ is determined by a Poisson random variable with parameter equal to the area of the region, and the number of points in non-overlapping regions are independent), and then defines $X$ as the solution to
\begin{align}\label{rep:2_homogeneous}
	X(t) = X(0) + \sum_{k=1}^K \zeta_k \int_{[0,t]\times [0,\infty)}1_{[q_{k-1}(s-),q_k(s-))}(x) N(ds\times dx)
\end{align}
where $q_k(t)=\sum_{\ell=1}^k \lambda_\ell(X(t)) $ with $q_0(t) \equiv 0$, and $1_{[a,b)}(x)$ is the indicator function that $x\in [a,b)$. See Figure \ref{fig:Space-TimeRepresentation_homogeneous} for a graphical realization of  representation  \eqref{rep:2_homogeneous} given three reaction channels. Simulation of this representation is equivalent to the well known Gillespie algorithm \cite{gillespie1976general,gillespie1977exact}.  Notice that this representation is equivalent to the stochastic differential equation form in \cite{li2007analysis}.

\begin{figure}
\begin{center}
 \includegraphics[width=3.8in]{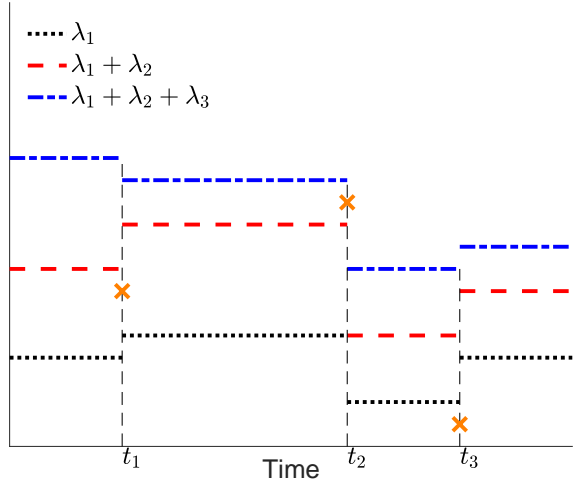} 
 \end{center}
 \caption{Graphical illustration of the space-time representation \eqref{rep:2_homogeneous}.  The $x$-axis represents time,  the orange marks denote the relevant points of the Poisson point process $N$.  
 The times at which the transitions occurred are labeled $t_1,t_2,$ and $t_3$.  
  The first transition occurred due to reaction 2 because the first point occurred under $\lambda_1(t_1) + \lambda_2(t_1)$, but above $\lambda_1(t_1)$.  The second transition was due to reaction 3, and the third transition was due to reaction 1.  Note that the intensity functions change only when reactions occur.}
 \label{fig:Space-TimeRepresentation_homogeneous}
\end{figure}

\item A third representation starts with independent uniform$[0,1]$ random variables $\{\xi_0,\xi_1,\dots\}$ and a unit-rate Poisson process $Y$, which is independent from $\{\xi_0,\xi_1,\dots\}$.  The process $X$ is then defined as the solution to 
 \begin{align}
 \begin{split}
 X(t) &=  X(0) + \sum_{k=1}^K \zeta_k \int_0^t  1_{   \left[ \frac{q_{k-1}(s-)}{\lambda_0(X(s-))}, \frac{q_{k}(s-) }{\lambda_0(X(s-))} \right)}( \xi_{R_0(s-)})dR_0(s)  \\
R_0(s) &=  Y\big( \int_0^t  \lambda_0(X(s)) ds\big), 
\end{split}
\label{rep:3_homogeneous}
\end{align}
where $\lambda_0(x) = \sum_{k=1}^K \lambda_k(x)$ and the $q_k$ are as above.  Note that simulation of the representation \eqref{rep:3_homogeneous} is also equivalent to Gillespie's algorithm \cite{gillespie1976general,gillespie1977exact},  though the processes \eqref{rep:2_homogeneous} and \eqref{rep:3_homogeneous} are built on different probability spaces. 
\end{enumerate}

 \subsection{The non-homogeneous case}

We turn to the non-homogeneous case.  The  modeling assumptions  analogous  to \eqref{eq:assump_jump_rates_homogeneous}  are
\begin{align}
\begin{split}
	P(\text{reaction $k$ occurs in $[t,t+h)$}\ | \ \mathcal{F}_t) &= \lambda_k(t,X(t)) h + o(h)\\
	P(\text{two or more reactions occur in $[t,t+h)$}\ | \ \mathcal{F}_t) &= o(h)\\
	P(\text{no reaction occurs in $[t,t+h)$}\ | \ \mathcal{F}_t) &=1 - \sum_{k}\lambda_k(t,X(t)) h + o(h),
	\end{split}
	\label{eq:assump_jump_rates_nonhomogeneous}
\end{align}
where the propensity functions are now functions of both time and the state of the system.

The  stochastic equations analogous to \eqref{eq:RTC_main}, \eqref{rep:2_homogeneous}, and \eqref{rep:3_homogeneous} are the following.
\begin{enumerate}
\item The   random time change representation  \eqref{eq:RTC_main} becomes
\begin{equation}\label{rep:1_nonhomogeneous}
X(t) = X(0)  + \sum_{k=1}^K Y_k \big( \int_0^t \lambda_k(s,X(s)) ds\big) \zeta_k,
\end{equation}
where the $Y_k$ are independent unit rate poisson processes.  

\item The space-time representation utilizing a Poisson point  process \eqref{rep:2_homogeneous} becomes 
\begin{align}\label{rep:2_nonhomogeneous}
X(t) &= X(0) + \sum_{k=1}^K \zeta_k \int_{[0,t]\times [0,\infty)}1_{[q_{k-1}(s-),q_k(s-))}(x) N(ds\times dx),
\end{align}
where $N$ is a space-time unit-rate Poisson point process and now $q_k(t)=\sum_{\ell=1}^k \lambda_\ell(t,X_t) $ for $k=1,2,....,K$ and $q_0(t) \equiv 0$. See Figure \ref{fig:Space-TimeRepresentation_nonhomogeneous} for a graphical realization of this representation given three reaction channels that is analogous to Figure \ref{fig:Space-TimeRepresentation_homogeneous}. 

\begin{figure}
\begin{center}
 \includegraphics[width=4.0in]{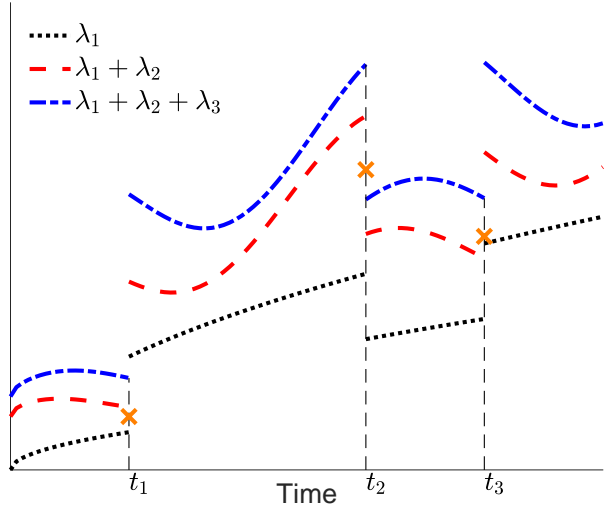} 
 \end{center}
 \caption{
  Graphical illustration of the space-time representation \eqref{rep:2_nonhomogeneous}.  The $x$-axis represents time,  the $\times$ marks denote the relevant points of the Poisson point process $N$.  
  The first transition occurred due to reaction 2 because the first point occurred under $\lambda_1(t_1) + \lambda_2(t_1)$, but above $\lambda_1(t_1)$. The second transition was due to reaction 2 as well, and the third due to reaction 3. }
 \label{fig:Space-TimeRepresentation_nonhomogeneous}
\end{figure}

\item The representation \eqref{rep:3_homogeneous} also stays largely the same.  We again start with independent uniform$[0,1]$ random variables $\{\xi_0,\xi_1,\dots\}$ and a unit-rate Poisson process $Y$, which is independent from $\{\xi_0,\xi_1,\dots\}$.  Then we define $X$ as the solution to
 \begin{align}
 \begin{split}
  X(t) &=  X(0) + \sum_{k=1}^K \zeta_k \int_0^t  1_{   \left[ \frac{q_{k-1}(s-)}{\lambda_0(s-,X(s-))}, \frac{q_{k}(s-) }{\lambda_0(s-,X(s-))} \right)}( \xi_{R_0(s-)})dR_0(s)  \\
R_0(s) &=  Y\big( \int_0^t  \lambda_0(s,X(s)) ds\big),
\end{split}
\label{rep:3_nonhomogeneous}
\end{align}
where $\lambda_0(s,x) = \sum_{k=1}^K \lambda_k(s,x)$ and the $q_k$ are as defined around \eqref{rep:2_nonhomogeneous}.  
\end{enumerate}

\subsection{Simulation and thinning}

Simulation of  each of the representations \eqref{rep:1_nonhomogeneous}, \eqref{rep:2_nonhomogeneous}, and \eqref{rep:3_nonhomogeneous}   consists of two steps:
\begin{enumerate}[(1)]
\item determine when the next transition will occur, and
\item determine which transition occurs at that time.
\end{enumerate}
To complete the two-step process, each of the representations requires the solution of a hitting time problem of the form \eqref{eq:int_to_estimate1}.  For example, in order to determine the time of the next transition in the representations \eqref{rep:2_nonhomogeneous} and \eqref{rep:3_nonhomogeneous}, we are required  to solve
\begin{equation}\label{eq:5678}
	\int_{t_\ell}^{t_\ell + \Delta} \lambda_0(s,X(s) ) ds = \mathcal E_{\ell},
\end{equation}
where $t_\ell$ is the time of the previous transition and $\mathcal E_{\ell}$ is a unit exponential random variable.  The representation \eqref{rep:1_nonhomogeneous} requires us to solve $K$ such equations and then take the minimum to determine which reaction causes the transition.  Numerically solving these hitting time problems is often numerically expensive as it requires a significant number of calculations to be made between each transition.

A well known method to avoid the numerical estimation of these integrals is to utilize thinning  \cite{GlynnAsmussen2007}.  Such methods have been introduced in the present context of biochemical processes in  \cite{lemaire2016exact, thanh2015simulation,voliotis2016stochastic}.  Here we show how these methods arise naturally via the representations presented above.  We will then be able to naturally develop  coupling methods for perturbed processes, which are the main contributions of this paper.


One natural way to think about thinning in the present context is to introduce an extra reaction channel, labeled as the $(K+1)$st, whose reaction vector is the zero vector, $\zeta_{K+1} = \vec 0\in \Z^d$.  Next,  choose $\overline \lambda_0$ satisfying 
\[
	\overline \lambda_0 \ge \lambda_0(s,X(s))=  \sum_{k=1}^K \lambda_k (s,X(s)),
\]
for all $s$ in some time-frame of interest.  Finally, take the following as the propensity for the $(K+1)$st reaction, 
\[
	\lambda_{K+1}(s,X(s)) = \overline \lambda_0 -\lambda_0(s,X(s)),
\]
for those $s$ in the time-frame of interest.  

Involving such a ``phantom reaction'' can be helpful because now   the previously time consuming process detailed around \eqref{eq:5678} for finding the time of the next reaction, $\Delta$, simply requires us to solve 
\[
	\int_{t_\ell}^{t_\ell+\Delta} \overline \lambda_0\, ds = \mathcal E_\ell \iff \Delta = \frac{\mathcal E_\ell}{\overline \lambda_0}.
\] 
That is, the determination of the time of the next reaction is exactly the same as for the standard Gillespie algorithm.  An extra cost comes, however,  in step (2) of the algorithm, which determines the reaction channel that causes a transition.  It  is now possible that the $(K+1)$st reaction will be chosen, in which case there is no update to the state of the system (though time is still flowed forward by $\Delta$ units).  In order to minimize the number of such ``wasted steps,'' it is desirable to select $\overline \lambda_0$ as close to $\sup_{s \in I} \lambda_0(s,X(s))$ as possible, where $I$ is our time-frame of interest.

Formally, representation \eqref{rep:2_nonhomogeneous} becomes 
\begin{align}\label{rep:2_nonhomogeneous_thin}
X(t) &= X(0) + \sum_{k=1}^{K+1} \zeta_k \int_{[0,t]\times [0,\infty)}1_{[q_{k-1}(s-),q_k(s-))}(x) N(ds\times dx),
\end{align}
 and representation \eqref{rep:3_nonhomogeneous} becomes
 \begin{align}\label{rep:3_nonhomogeneous_thin}
 \begin{split}
  X(t) &=  X(0) + \sum_{k=1}^{K+1} \zeta_k \int_0^t  1_{   \left[ \frac{q_{k-1}(s-)}{\overline \lambda_0}, \frac{q_{k}(s-) }{\overline \lambda_0} \right)}( \xi_{R_0(s-)})dR_0(s)  \\
R_0(s) &=  Y\big( \int_0^t  \overline \lambda_0 ds\big),
\end{split}
\end{align}
with all notation detailed above.

An algorithm for the numerical  implementation of either representation is the following.  The presented algorithm is essentially  the same as the \textit{Extrande} method in \cite{voliotis2016stochastic}.  Note that the algorithm below is slightly different than the representations above in that we update $\overline \lambda_0$ after each step.  This is left out of the representations for notational clarity.
 
\begin{algorithm}[Simulation of \eqref{rep:2_nonhomogeneous_thin} or \eqref{rep:3_nonhomogeneous_thin}]\label{algorithm1}
Initialize by setting $t=0$ and $X = X(0)$.  Let $T>0$ be some predetermined end time.  Set $q_0 \equiv 0$.   Do the following until $t$ exceeds  $T$.    All generated random variables are independent of previously generated random variables.
\begin{enumerate}
\item  Let $\overline \lambda_0$ satisfy $\overline \lambda_0 \ge\sum_{k=1}^K \lambda_k(s,X)$ for all $s \in [t,T]$.  
\item Generate a unit exponential random variable $\mathcal E$ and set $\Delta = \frac{\mathcal E}{\overline \lambda_0}.$
\item  If $t +\Delta> T$, do the following
\begin{enumerate}[(i)]
\item Define $X(u) = X$ for $u \in [t, T]$,
\item set $t = T$ and end loop.
\end{enumerate}
Otherwise do the following:  

\item  Calculate the propensities $\lambda_k(t+\Delta,X)$ and $q_k(t+\Delta) = \sum_{\ell=1}^k \lambda_l(t+\Delta,X)$, for $k \in \{1,\dots,K+1\}$.
\item \textit{[Thinning step]} Generate a Uniform$[0,1]$ random variable $\xi$, and find  $\mu\in \{1, 2, \dots, K, K+1\}$ for which  
\[
	\frac{q_{\mu-1}(t+\Delta)}{\overline \lambda_0}\le   \xi <  \frac{q_{\mu}(t+\Delta) }{\overline \lambda_0}.
\]
 
\item Set 
\begin{enumerate}[(i)]
\item $X(u) = X$ for $u \in [t, t+\Delta)$,
\item  $X \leftarrow X + \zeta_\mu$, 
\item $t\leftarrow t+\Delta$,
\end{enumerate}
and go to step 1. 
\end{enumerate}
\end{algorithm}

\begin{figure}[b!]
     \centering
      \includegraphics[width=3.88in]{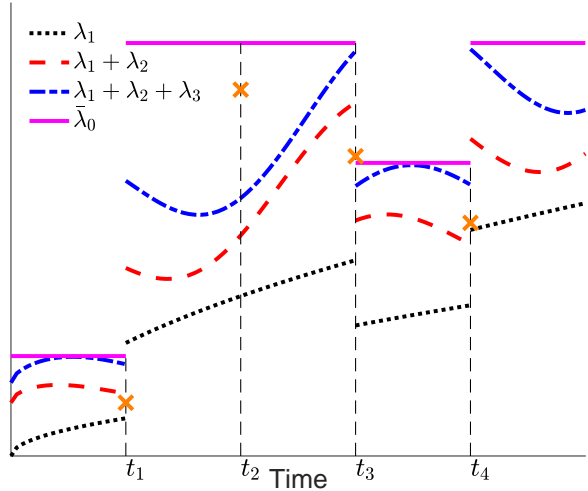} 
      \caption{Graphical illustration of the Extrande algorithm.     The $x$-axis represents time,  the $\times$ marks denote the relevant points of the Poisson point process $N$.  
  The first transition occurred due to reaction 2 since the first point occurred under $\lambda_1(t_1) + \lambda_2(t_1)$, but above $\lambda_1(t_1)$. There is no transition at the second time $t_2$ as the point was between $\lambda_1(t_2) + \lambda_2(t_2) + \lambda_3(t_2)$ and $\overline \lambda_0$.  The transitions at times $t_3$ and $t_4$ were due reactions 2 and 3, respectively.}
      \label{Extrande_Graph}
\end{figure}

\begin{remark}\label{remark1}
 It could be the case that   $ \sup_{s\in [t,T]}\sum_{k=1}^K \lambda_k(s,X)$ is much larger than $\sum_{k=1}^K \lambda_k(t,X)$.  When this occurs,  Algorithm \ref{algorithm1} will be  inefficient because the probability of acceptance in step 5 will be low (i.e., we will have $\mu = K+1$ and $\zeta_\mu = \vec 0$ with large probability). 
 When this occurs we may simply choose a $\overline \lambda_0$ to be some value greater than $\lambda_0(t,X)$ and then stop the step  if the intensity goes above $\overline \lambda_0$.  
 To be specific, we would let $t_0=  \inf\{ u: \sum_{k=1}^K \lambda_k(u,X)  \ge \overline \lambda_0\}$ and if $\Delta > t_0 - t$ we will update the process by setting $t = t_0$ and ``update'' the state vector by adding $\zeta_{K+1} = \vec 0$.  We would then return to step 1 of the algorithm.
\end{remark}
A graphical realization of the algorithm as applied to the representation \eqref{rep:2_nonhomogeneous_thin} in which $\overline \lambda_0$ is updated after each step, and which is analogous to Figures \ref{fig:Space-TimeRepresentation_homogeneous} and \ref{fig:Space-TimeRepresentation_nonhomogeneous},  can be found in  Figure \ref{Extrande_Graph}. 


 \section{Couplings}
 \label{sec:couplings1}

To motivate the idea of coupling two processes, we introduce one example concerning parametric sensitivities, or derivatives of expectations with respect to parameters. Since population processes can be modeled as reaction networks, consider the following, which arises from an ecological model in \cite{peel2014effect} where it was used to investigate the persistence of a pathogen in wild mammalian populations.

\begin{model}\label{example:SIR}
Consider an SIR model with the following possible transitions: 
\begin{align*}
(S,I,R) &\overset{\lambda_1(t)}{\to} (S+1,I,R) \tag{R1} \\
S &\overset{m}{\to} \emptyset  \tag{R2} \\
I &\overset{m}{\to} \emptyset   \tag{R3} \\
R &\overset{m}{\to} \emptyset \tag{R4} \\
S+I &\overset{\beta}{\to}  2I   \tag{R5} \\
I &\overset{\gamma}{\to} R. \tag{R6}
\end{align*}
Each of the reactions except for $(R1)$ and $(R5)$ is assumed to have mass-action kinetics.  For $(R1)$ it is assumed that  $\lambda_1(t) = B(t)(S(t)+I(t)+R(t))$ where $B(t)$ is a periodic  function that  models  instantaneous birth pulses in a limited period of the year. It takes the form:
\begin{equation*}
B(t) = k \exp{(-s \cos(\pi t-\phi)^2)},
\end{equation*}
where $s$ is the synchrony parameter that controls the duration of the birth pulse and $\phi$ is the phase parameter which determines the timing of birth peaks. $k$ is a scaling constant which is determined by $k = \frac{me^{s/2}}{I_0(s/2)}$, where $m$ is the death rate (which is also used in reactions $(R2),$ $(R3)$, and $(R4)$) and $I_0(z) = \frac{1}{\pi} \int_0^\pi e^{z\cos{x}}dx$ is the modified 
Bessel equation of the first kind. 
For $(R5)$, $\beta$ is the transmission rate and the propensity/rate of the reaction is $\lambda_5(t) = \frac{\beta S(t)I(t)}{(S(t)+I(t)+R(t))}$.  Finally, note that  $\gamma$ is the recovery rate of the infected individuals. \newline

In this model, the quantity we are interested in is the sensitivity of the extinction probability of pathogens within 10 years against all the parameters in our model. More specifically, we define the stopping time $\displaystyle \tau = \inf_{t\ge 0} \{ I(t) = 0\}$ and the extinction probability $p = E[1_{\tau<10}] =  P(\tau<10)$. The quantities we are looking for are the derivatives of $p$ against all sets of parameters $m$, $\gamma$, $\phi$, $R_0$ and $s$, where $R_0$ is the basic reproduction ratio satisfying $\beta = R_0(m+\gamma)$.  In particular, we would be interested in determining which parameters are the dominant factors in extinction events.   \hfill $\square$
\end{model}

One common method to calculate parametric sensitivities is to use a finite difference method.  Let  $\theta$ be some parameter of interest in the model, and let the parameterized process be denoted by $X^\theta$.  Suppose that $E[f(X^\theta)]$ is some expectation of interest, where $f$ is some path functional.  Then utilize the approximation
\begin{align*}
\frac{d}{d\theta} E[f(X^\theta)] \approx \frac{E[f(X^{\theta+h/2}) ]- E[f(X^{\theta-h/2})]}{h} = E\bigg[\frac{f(X^{\theta+h/2})-f(X^{\theta-h/2})}{h}\bigg],
\end{align*}
which is a centered finite difference.
The natural Monte Carlo estimator is then
\[
	\frac1n \sum_{i = 1}^n \frac1h (f(X^{\theta+h/2}_{[i]}) - f(X^{\theta-h/2}_{[i]})),
\]
where $(X^{\theta+h/2}_{[i]}, X^{\theta-h/2}_{[i]})$ is the $i$th independent pair of processes generated.  As is always the case with Monte Carlo methods, the efficiency of the method scales directly with the relevant variance, which in this case is
\[
	\text{Var}\left(f(X^{\theta+h/2}_{[i]}) - f(X^{\theta-h/2}_{[i]})\right).
\]

We leave the setting of finite differences and return to simply considering two processes, $X$ and $Z$, which share the reaction vectors $\zeta_k$, but differ in their propensity functions. (For example, in the discussion above $X$ could be $X^{\theta+h/2}$ and $Z$ could be $X^{\theta-h/2}$.) We will denote the propensity functions of $X$ by  $\lambda_k^X (s,X(s))$ and the propensity functions of $Z$ by $\lambda_k^Z (s,Z(s))$.  We will denote the sums of the propensity functions by 
\[
	\lambda_0^X(s,X(s)) = \sum_{k = 1}^K \lambda_k^X(s,X(s))\quad \text{ and } \lambda_0^Z(s,Z(s)) = \sum_{k = 1}^K \lambda_k^Z(s,Z(s)).
	\]
In this section, we will introduce four different couplings of the jump processes $(X,Z)$.

\subsubsection*{Coupling \#1: Independent samples} The easiest way to get the correct marginal processes is to simply generate independent realizations of $X$ and $Z$.  For example, this could be done by implementing Algorithm \ref{algorithm1} for the two process with independent uniform and exponential random variables.  Since the processes are \textit{uncoupled}, we expect any other coupling to provide  lower variances for the requisite Monte Carlo estimators.

\subsubsection*{Coupling \#2: ExtrandeCRN} A natural coupling, which we will term \textit{ExtrandeCRN}, is to utilize the same seed of your random number generator for the implementation of Algorithm \ref{algorithm1} for both $X$ and $Z$.  This method is often termed ``common random numbers.'' Mathematically, this coupling is equivalent to using the same Poisson process $Y$ and uniform random variables $\{\xi_0,\xi_1,\dots\}$ in the representation \eqref{rep:3_nonhomogeneous_thin}.

Note that for this coupling the processes are not guaranteed to have the same jump times, nor are they guaranteed to have the same sequence of transitions.

\subsubsection*{Coupling \#3: Extrande thinning}  
\begin{figure}
     \centering
      \includegraphics[width=5in]{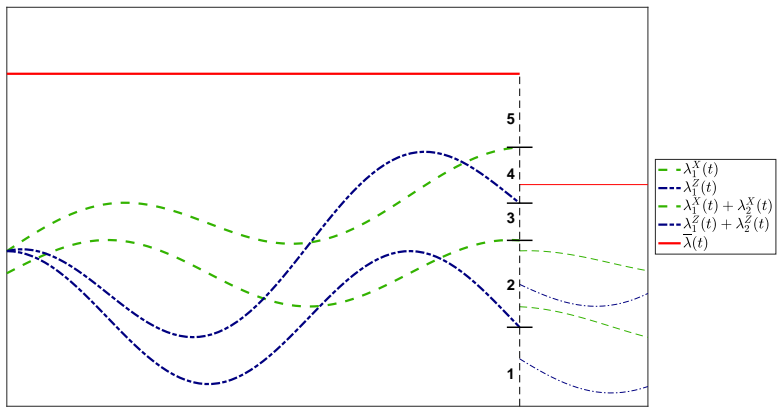} 
      \caption{Graphical illustration of the \textit{Extrande thinning} coupling for processes $X$ and $Z$ with  two reaction channels.  The dashed green curves determine the propensity functions for $X$, whereas the dash-dotted curves determine the propensities for $Z$.    The solid red line is $\overline \lambda_0$.  Suppose a point of the point process $N$ occurs along the vertical dashed line below $\overline \lambda_0$.  There are five possibilities.  If the point occurs in region 1, both  $X$ and $Z$  will transition via reaction 1.  If the point occurs in region 2,  $X$  will transition via reaction 1 and $Z$ will transition via reaction 2.   If the point occurs in region 3, both  $X$ and $Z$  will transition via reaction 2. If the point occurs in region 4,  $X$  will transition via reaction 2 and  $Z$  will not transition.  If the point occurs in region 5, neither $X$ nor $Z$ will transition.}
      \label{Extrande_thinning}
\end{figure}

The previous coupling arose by generating $X$ and $Z$ via \eqref{rep:3_nonhomogeneous_thin} with the same sources of randomness (the Poisson process, via the unit exponentials, and the uniforms).  Our third coupling, termed \textit{Extrande thinning}, arises by generating $X$ and $Z$ via \eqref{rep:2_nonhomogeneous_thin} with the same sources of randomness.  Specifically, $X$ and $Z$ will each satisfy   \eqref{rep:2_nonhomogeneous_thin} with  the \textit{same} unit-rate Poisson point process, $N$.   Figure \ref{Extrande_thinning} offers a graphical illustration.



\begin{algorithm}[Extrande thinning]
Initialize by setting $t=0$, $X = X(0)$, and $Z=Z(0)$.  Let $T>0$ be some predetermined end time.  Do the following until $t$ exceeds $T$.  All generated random variables are independent of previously generated random variables.

\begin{enumerate}[1.]
\item Let  $\overline \lambda_0$ satisfy $\overline \lambda_0 \ge \sup_{s \in [t,T]} \{\lambda_0^X(s,X), \lambda_0^Z(s,Z)\}$.
 
\item Generate a unit exponential random variable $\mathcal E$ and set $\Delta = \frac{\mathcal E}{\overline \lambda_0}.$

\item If $t + \Delta > T$, do the following
\begin{enumerate}[(i)]
\item Define $X(u) = X$ and $Z(u) = Z$ for $u \in [t,T]$,
\item set $t = T$ and end loop.
\end{enumerate}
Otherwise do the following:
\item  Calculate the propensities $\lambda_k^X(t+\Delta,X)$, $\lambda_k^Z(t+\Delta,Z)$ and  for $k\in \{1,\dots,K+1\}$,
\[
	q_k^X(t+\Delta) =\sum_{\ell=1}^k \lambda_\ell^X(t+\Delta,X) \quad \text{ and } \quad q_k^Z(t+\Delta)= \sum_{\ell=1}^k \lambda_\ell^Z(t+\Delta,Z),
\]
and $q_0^X(t+\Delta) = q_0^Z(t+\Delta) = 0$.
\item \textit{[Thinning step]} Generate a Uniform$[0,1]$ random variable $\xi$, find  $\mu, \nu\in \{1,\dots,K+1\}$ such that  
\[
\frac{q_{\mu-1}^X(t)}{\overline \lambda_0} \le \xi < \frac{q_{\mu}^X(t) }{\overline \lambda_0}  \quad \text{and}\quad \frac{q_{\nu-1}^Z(t)}{\overline \lambda_0} \le \xi < \frac{q_{\nu}^Z(t) }{\overline \lambda_0}.
\]
(Recall that $\zeta_{K+1} = \vec 0$ for each process.)
\item Set
\begin{enumerate}
\item $X(u) = X$ and $Z(u) = Z$ for $u\in [t,t+\Delta)$,
\item $X \leftarrow X + \zeta_\mu$ and $Z \leftarrow Z + \zeta_\nu$,
\item $t \leftarrow t + \Delta$.
\end{enumerate}
and go to step 1.
\end{enumerate}
\label{algorithm2}
\end{algorithm}

\begin{remark}\label{remark2}
Just as was pointed out in Remark \ref{remark1} pertaining to Algorithm \ref{algorithm1}, if there is no $\overline \lambda_0$ satisfying the inequality in step 1 that is close in value to $\max\{\lambda_0^X(t,X),\lambda_0^Z(t,Z)\}$, then you may choose $\overline \lambda_0$ to be some value greater than $\max\{\lambda_0^X(t,X),\lambda_0^Z(t,Z)\}$ and simply stop the step at the time one of the total intensities goes above $\overline \lambda_0$. 
\end{remark}

\subsubsection*{Coupling \#4: The stacked coupling} 

In ExtrandeCRN (Coupling \#2), the coupled processes will not, in general, have transitions at the same time.  On the other hand, Extrande thinning (Coupling \#3) will often yield simultaneous transitions for  $X$ and $Z$.  However, it can be the case that even when $X$ and $Z$ transition at the same time in Extrande thinning, the reaction channels responsible for those transitions will be different  (see region 2 in Figure \ref{Extrande_thinning}).  Our final coupling, the \textit{stacked coupling}, improves the situation again.  For this coupling the processes will often have transitions at the same time, and when both processes transition simultaneously they necessarily transition according the same reaction channel.

\begin{figure}
     \centering
      \includegraphics[width=5.5in]{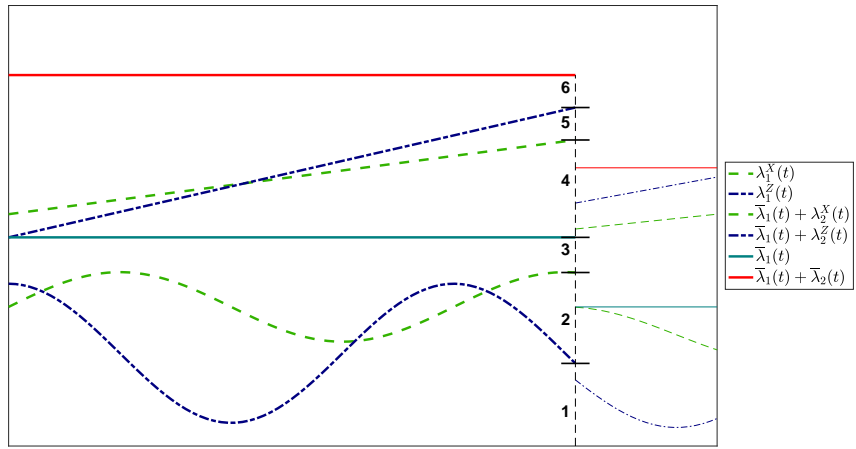} 
      \caption{Graphical illustration of the  \textit{stacked coupling} for processes $X$ and $Z$ with  two reaction channels.  The dashed green and dash-dotted blue curves are related to the propensity functions for $X$ and $Z$, respectively.  For $k \in \{1,2\}$,  $\overline \lambda_k$  gives an upper bound for $\sup_s\{\lambda_k^X(s,X(s)),\lambda_k^Z(s,Z(s))\}$ over our time period of interest.  The solid red line on the top is $\overline \lambda_0^{\text{total}} = \overline \lambda_1 + \overline \lambda_2$.  Suppose a point of the point process  occurs along the vertical dashed line below $\overline \lambda_0$.  There are six possibilities.  If the point occurs in region 1, both  $X$ and $Z$  will transition via reaction 1.  If the point occurs in region 2,  $X$  will transition via reaction 1 and $Z$ will not transition.   If the point occurs in region 3, neither $X$ nor $Z$ will transition. If the point occurs in region 4,  both $X$  and $Z$ will transition via reaction 2.  If the point occurs in region 5, $Z$ will transition via reaction 2 and $X$ will not transition.  If the point occurs in region 6, neither $X$ nor $Z$ will transition.}
      \label{Stacked_coupling}
\end{figure}

  The stacked coupling comes from giving each reaction channel its own piece of $\R^2_{\ge 0}$ in the space-time representation.    It is essentially a space-time, and thinned, version of the coupled finite difference method introduced in \cite{AndCFD2012} in that the propensity for a particular reaction channel  to cause  simultaneous transitions for  $X$ and $Z$ is the minimum of the respective propensities.  See Figure \ref{Stacked_coupling}.

The stochastic representation of this coupling can be written in this form: 
\begin{align}
\label{eq:stacked_rep}
\begin{split}
X(t) = X(0) + \sum_{k=1}^K  \zeta_k & \int_{[0,t]\times[0,\infty)} 1_{\left[ \overline q_{k-1},\ \overline q_{k-1}+ \lambda_k^X(s-,X(s-) \right)}(x) N(ds\times dx) \\
Z(t)= Z(0) + \sum_{k=1}^K  \zeta_k & \int_{[0,t]\times[0,\infty)} 1_{\left[ \overline q_{k-1},\ \overline q_{k-1}+ \lambda_k^Z(s-,Z(s-) \right)}(x) N(ds\times dx), 
\end{split}
\end{align}
where
\begin{align*}
	\overline q_k = \sum_{\ell=1}^k \overline \lambda_\ell, \quad \text{ with } \quad \overline \lambda_\ell \ge  \sup_{s\in I}\{\lambda_l^X(s,X(s)),\lambda_l^Z(s,Z(s))\},
\end{align*}
where $I$ is our time period of interest and $\overline q_0 = 0$.  As a final bit of notation, we let
\[
	 \overline \lambda_0^{\text{total}} =  \sum_{\ell =1}^K \overline \lambda_\ell.
\]

%

Note that the algorithm presented below is slightly different than the representation \eqref{eq:stacked_rep} as in the algorithm we update $\overline \lambda_k, \overline q_k,$ and $\overline \lambda_0^{\text{total}}$ after each step.  In \eqref{eq:stacked_rep} these terms are fixed.

\begin{algorithm}[Stacked coupling]
Initialize by setting $t=0$, $X = X(0)$, and $Z=Z(0)$.  Let $T>0$ be some predetermined end time. Set $\overline q_0\equiv 0$.  Do the following until $t$ exceeds $T$.  All generated random variables are independent of previously generated random variables.

\begin{enumerate}[1.]
\item Find all $\overline \lambda_k$, $\overline q_k$, and $\overline \lambda_0^{\text{total}}$ where $[t,T]$ is the time period of interest.

 
\item Generate a unit exponential random variable $\mathcal E$ and set $\Delta = \frac{\mathcal E}{\overline \lambda_0^{\text{total}}}.$

\item If $t + \Delta > T$, do the following
\begin{enumerate}[(i)]
\item Define $X(u) = X$ and $Z(u) = Z$ for $u \in [t,T]$,
\item set $t = T$ and end loop.
\end{enumerate}
Otherwise do the following:
\item Set $X(u) = X$ and $Z(u) = Z$ for $u \in [t,t+\Delta)$.
\item  Calculate the propensities $\lambda_k^X(t+\Delta,X)$, $\lambda_k^Z(t+\Delta,Z)$  for $k\in \{1,\dots,K\}$.

\item \textit{[Thinning step]} Generate a uniform random variable $\xi$ in $[0,1]$.
\begin{enumerate}[(i)]
\item Find  $\mu \in \{1,\dots,K\}$ for which
\[
\frac{\overline q_{\mu-1}}{\overline \lambda_0^{\text{total}}} \le  \xi < \frac{ \overline q_{\mu} }{\overline \lambda_0^{\text{total}}}.
\]

\item  If 
\[
\frac{\overline q_{\mu-1}}{\overline \lambda_0^\text{total}} \le  \xi < \frac{\overline q_{\mu-1} +  \lambda_\mu^X(t+\Delta,X)}{\overline \lambda_0^{\text{total}}},
\]
 then set $X \leftarrow X + \zeta_\mu$.
\item If 
\[
\frac{\overline q_{\mu-1}}{\overline \lambda_0^\text{total}} \le  \xi < \frac{\overline q_{\mu-1} +  \lambda_\mu^Z(t+\Delta,Z)}{\overline \lambda_0^{\text{total}}},
\]
 then set $Z \leftarrow Z + \zeta_\mu$.
\end{enumerate}

\item Set $t \leftarrow t + \Delta$ and go to step 1. 
\end{enumerate}

\label{algorithm3}
\end{algorithm}

\begin{remark}\label{remark3}
Similar to Remarks \ref{remark1} and \ref{remark2}, the upper bounds in step 1 of the algorithm do not have to be uniform over the time period $[t,T]$ so long as you stop the step if one of the relevant intensities goes above the chosen upper bound.
\end{remark}

\section{Applications and examples}
\label{sec:applications}

Two application areas that benefit greatly from good coupling methods are parametric sensitivities and multilevel Monte Carlo.  We briefly introduced finite difference methods for parametric sensitivities in Example \ref{example:SIR} of Section \ref{sec:couplings1}.  We will briefly introduce multilevel Monte Carlo in the present setting, but point the reader to  \cite{AndHigham2012,AHS2014,lester2016extending} for a deeper introduction.  We will then introduce two more example models and then provide numerical tests for the devised coupling methods.

\subsection{Multilevel Monte Carlo in the chemical kinetic context}

Multilevel Monte Carlo (MLMC) is a method that uses sequences of coupled processes in order to efficiently compute expectations \cite{Giles2008}.    Let $f$ be some path functional of interest, let $X$ be a relevant stochastic process and consider the problem of numerically approximating the expectation $E[f(X)]$.

The MLMC estimator of \cite{AndHigham2012} is built in the following manner.   For a fixed integer $M$, and $\ell \in \{\ell_0,\ell_0 + 1,\dots,L\}$, where both $\ell_0$ and $L$ depend upon the model and path-wise simulation method being used, let  $\hl = T M^{-\ell}$, where $T$ is our terminal time.  Reasonable choices for $M$ include integers between 2 and 7.  Now note 
\begin{align}
  E [f(X)] =  E [f(X) - f(Z_{L})] +  \sum_{\ell = \ell_0 + 1}^L E[ f(Z_{\ell}) - f(Z_{\ell - 1})] + E [ f(Z_{\ell_0}) ], 
  \label{eq:add_sub}
\end{align}
where $Z_\ell$ is a realization generated via tau-leaping  \cite{Anderson2007b,Gill2001} (i.e.~Euler's method) with a time-step of $h_\ell$, and where the telescoping sum is the key feature to note.  We define independent estimators for the multiple terms by
\begin{align*}
	\Ql_{E}\ &= \ \frac{1}{n_{E}} \sum_{i = 1}^{n_{E}} (f(X_{[i]}) - f(Z_{L,[i]})) ,\ \ \ \Ql_{\ell} \ = \ \frac{1}{n_{\ell}} \sum_{i = 1}^{n_{\ell}} ( f(Z_{\ell,[i]}) - f(Z_{\ell - 1,[i]})),\\
		\Ql_{\ell_0}\ &= \ \frac{1}{n_{0}} \sum_{i = 1}^{n_{0}} f(Z_{\ell_0,[i]}), 
\end{align*}
where  $\ell \in \{\ell_0 + 1,\dots,L\}$ and $n_E,n_0,n_\ell$ must be determined during the computation.  Now note that 
\begin{equation}
	\Ql  =  \Ql_E +  \sum_{\ell = \ell_0}^L \Ql_{\ell}
	\label{eq:unbiased_MLMC}
\end{equation}
is an \textit{unbiased} estimator for $E[ f(X)]$.   The above observations are not useful unless we can successfully \textit{couple} the process $(X,Z_L)$ and $(Z_\ell,Z_{\ell-1})$  in a way that significantly reduces the variance of the estimator $\hat Q$ at each level.  That is, we want to minimize $\text{Var}(f(X) - f(Z_L))$ and $\text{Var}(f(Z_{\ell}) - f(Z_{\ell-1}))$.  

Even though we are considering non-homogeneous processes, the Euler approximations $Z_\ell$ have intensities that are constant over time-steps of size $h_\ell$ at $\lambda_k(t_n,Z_\ell(t_n))$, where we assume a current time of $t_n$.  Hence, the coupled processes $(Z_\ell, Z_{\ell-1})$ can be efficiently generated via  the standard coupling introduced in \cite{AndHigham2012}.  However,  the coupled processes $(X,Z_\ell)$ are non-homogeneous (because of $X$), and therefore we will utilize the stacked coupling.

\subsection{Example models}
We will test our couplings on three models.   Example Model \ref{example:SIR} of Section \ref{sec:couplings1} is one of the models.  We provide two more here.
\label{sec:example}

\begin{model}
\label{example:dimer}
We consider a  model of transcription-translation and dimerization that generalizes Example Model \ref{example:first}:
\begin{align*}
	\emptyset &\overset{\lambda_1(t)}{\to}  M  \tag{R1}\\
	M &\overset{100}{\to} M + P\tag{R2}\\
	M &\overset{1}{\to} \emptyset\tag{R3}\\
	P &\overset{1}{\to} \emptyset\tag{R4}\\
	2P &\overset{3\cdot 10^{-7}}\to D\label{eq:R5}\tag{R5}\\
	D &\overset{10}\to \emptyset,\label{eq:R6}\tag{R6}
\end{align*}
with
\[
	\lambda_1(t)=60+15\sin(\tfrac{2\pi t}{24}).
\]
We will take the following as our initial condition
\[
	X(0) = (0, 1000,0),
\]
with the species ordered as $M, P$, and $D$.
 \hfill $\square$
\end{model}

\begin{model}
\label{example:MMP}
Consider the network
\begin{align*}
S_1 + S_2 &\overset{\lambda_1(t)}{\to} S_3 \tag{R1} \\
S_3 &\overset{1}{\to} S_1 + S_2     \tag{R2} \\
S_3 &\overset{1}{\to} S_2 + S_4      \tag{R3}
\end{align*} 
with $\lambda_1(t)  = \frac{k_1(t)}{1000}$, and where $k_1(t)$  is the state of a Markov process with state space $\{ 0.5,1.5,5\}$,  unit exponential holding times, and transition probability matrix 
\[
\left[
\begin{array}{ccc}
0 & 0.5 & 0.5 \\
0.5 & 0 & 0.5 \\
0.5 & 0.5 & 0 \\ \end{array}\right].
\]
We will take as our initial condition,
\[
	X(0) = (1000, 1000, 0, 0)\quad \text{ and } \quad k_1(0) = 0.5.
	\]
	\hfill $\square$
\end{model}

\subsection{Numerical examples for parametric sensitivities}
\begin{example}\label{numexample1}
We consider Example Model \ref{example:dimer} with propensity function
\[
	\lambda_1^\theta(t) = 60 +  \theta \sin(\tfrac{2\pi t}{24})
\]
for reaction channel (R1).  We will parameterize the model by $\theta$, denoting the state of the system at time $t$ by $X^\theta(t)$.  In this example, we are first interested in estimating the sensitivity of the expected number of mRNA molecules with respect to $\theta$ for $t\in [0,20]$, at $\theta = 15$.  We will use centered finite differences
\begin{align}\label{eq:toestimate1}
	\frac{d}{d\theta}E[X_{1}^\theta(t)] \big|_{\theta= 15}\approx \frac1hE\big[ X_{1}^{15+h/2}(t) -X_{1}^{15-h/2}(t)  \big]  \quad\quad \text{         for  $t\in [0,20]$}
\end{align}
with the four couplings presented in Section \ref{sec:couplings1}.

To couple the process $X^{\theta+h/2}$ and $X^{\theta-h/2}$, it is useful to notice that $\lambda_1^{\theta+h/2}$ and $\lambda_1^{\theta-h/2}$ are bounded above by $60+\theta + \frac{h}{2}$ and $60+\theta - \frac{h}{2}$, respectively.   The other reaction channels have fixed parameters and their values (and hence upper bounds) are fixed between transitions.

To gather statistics, we generated 50,000 coupled paths with $h=0.1$ using the different methods presented in Section \ref{sec:couplings1} in order to estimate the right-hand side of \eqref{eq:toestimate1}. 
\begin{itemize}
\item  In Figure \ref{fig:senmod1}  we plot  the estimated sensitivities \eqref{eq:toestimate1} of the different estimators as functions of time.  We note that each captures the same overall behavior, though Extrande thinning and  the stacked coupling are significantly less variable.
\begin{figure}
     \centering
      \includegraphics[width=3in]{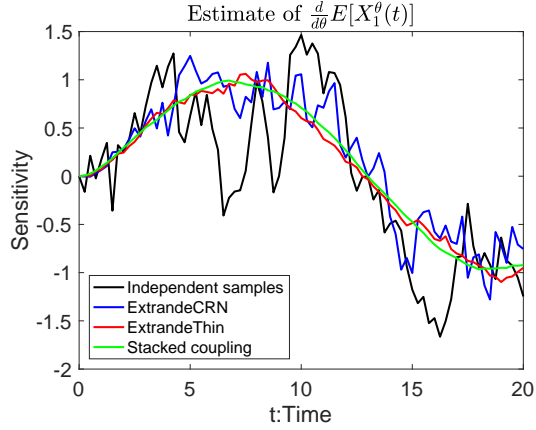} 
      \caption{Estimation of $\frac{d}{d\theta}E[X_{1}^\theta(t)]$ at $\theta = 15$ using finite differences and four different couplings for the model in Example \ref{numexample1}.}
      \label{fig:senmod1}
\end{figure}

\item In Figure \ref{fig:senmod2} we plot the estimate for $t\in [0,20]$:
\[
	\text{Var}\bigg(\frac{X_{1}^{\theta+h/2}(t) -X_{1}^{\theta-h/2}(t))}{t}\bigg)
\]
at $\theta = 15$ for each of the different methods.  We see that the variance of the stacked coupling is dramatically lower than the variance of the other couplings.

\begin{figure}
\begin{center}
\subfigure{\includegraphics[width = 3in]{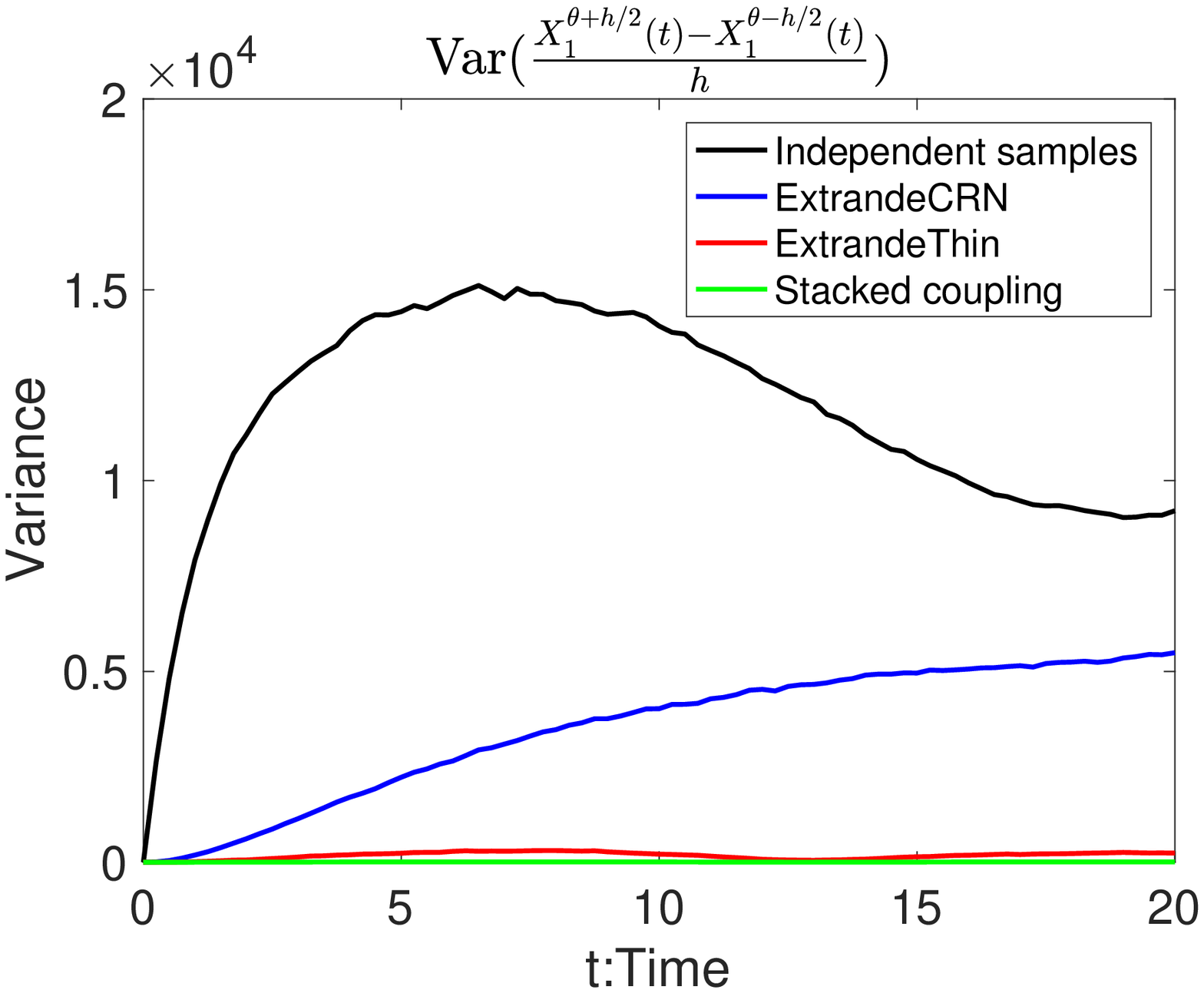}}
\subfigure{\includegraphics[width = 3in]{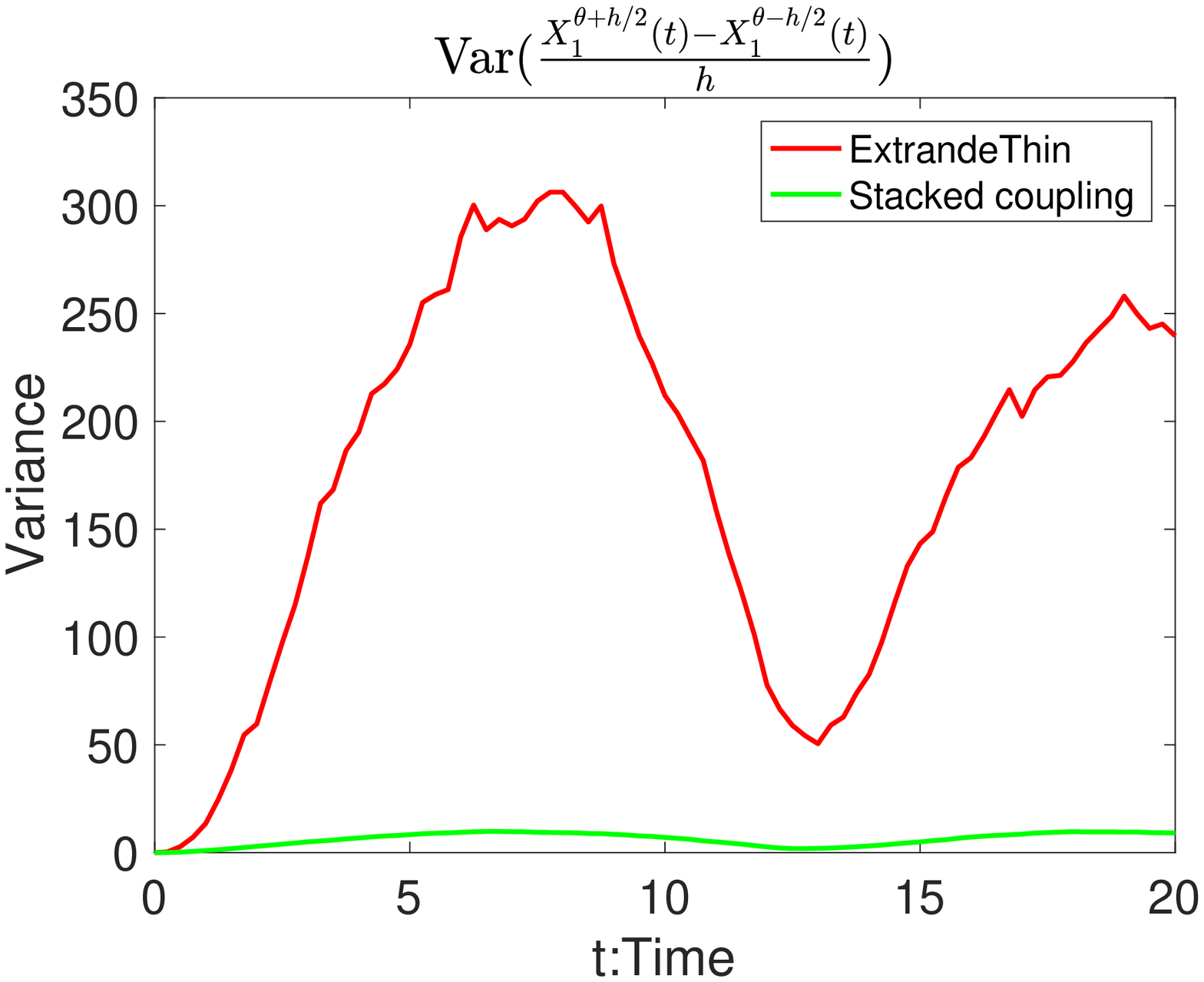}} 
\end{center}
\caption{Estimates of $\text{Var}\big(\frac{X_{1}^{\theta+h/2}(t) -X_{1}^{\theta-h/2}(t))}{h}\big)$ at $\theta = 15$ for the different couplings.}
\label{fig:senmod2}
\end{figure}
\end{itemize}

We now turn to estimating the sensitivity of the dimers:
\begin{align}\label{eq:toestimate_dimers}
	\frac{d}{d\theta}E[X_{3}^\theta(t)] \big|_{\theta= 15}\approx \frac1hE\big[ X_{3}^{15+h/2}(t) -X_{3}^{15-h/2}(t) \big] \quad\quad \text{         for  $t\in [0,20]$}
\end{align}
We again use a perturbation of $h=0.1$ and 50,000 coupled sample paths for each method.

\begin{itemize}
\item  In Figure  \ref{fig:senmod3} we plot  the estimated sensitivities of the different estimators as functions of time. 
\begin{figure}
     \centering
      \includegraphics[width=3in]{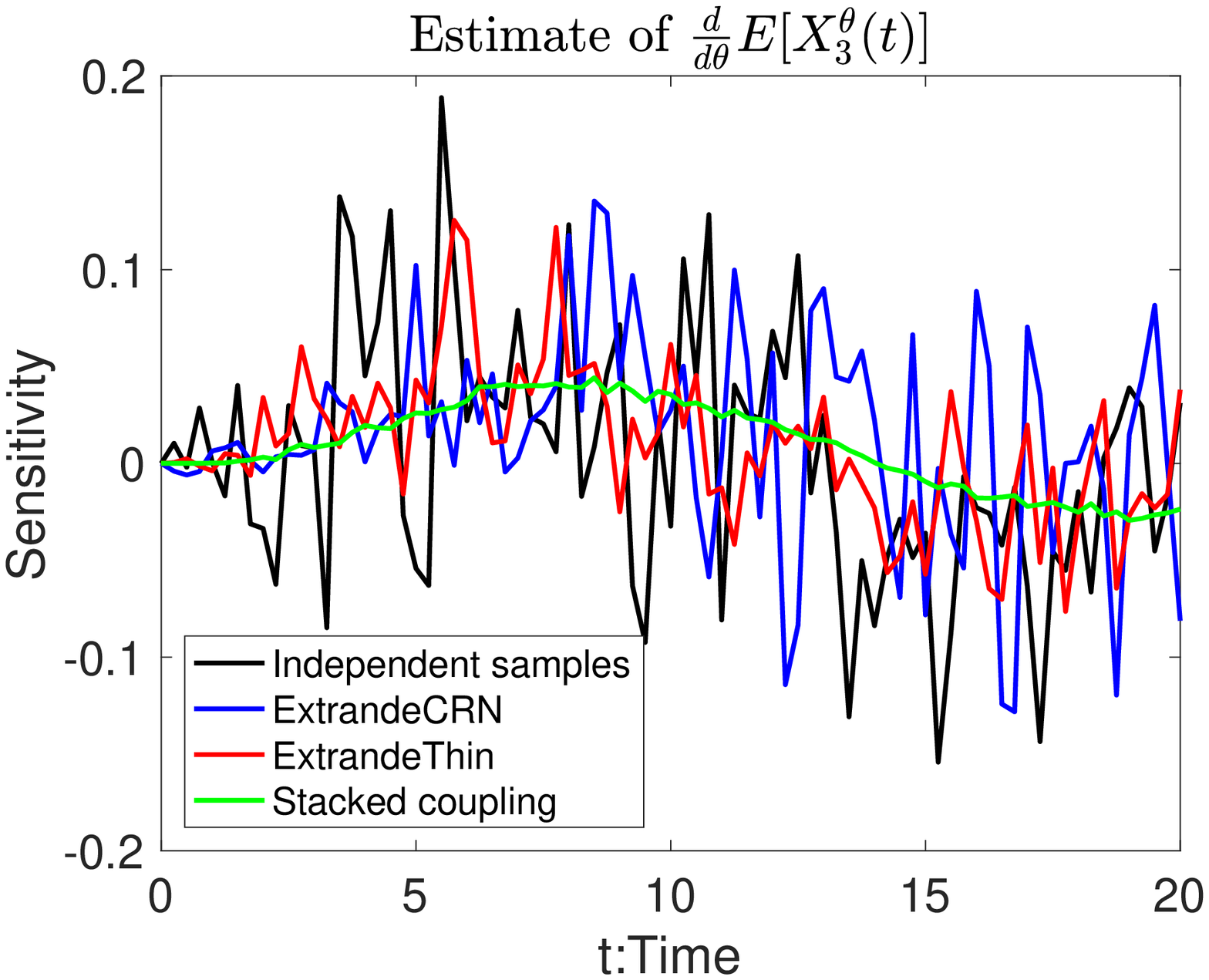} 
      \caption{Estimation of $\frac{d}{d\theta}E[X_{3}^\theta(t)]$ at $\theta = 15$ using finite differences and four different couplings  for the model in Example \ref{numexample1}.}
      \label{fig:senmod3}
\end{figure}

\item In Figure \ref{fig:senmod4} we plot the estimate for 
\[
	\text{Var}\bigg(\frac{X_{3}^{\theta+h/2}(t) -X_{3}^{\theta-h/2}(t))}{h}\bigg)
\]
at $\theta = 15$
for each of the different methods.  We again see that the variance of the stacked coupling is dramatically lower than the variance of the other couplings.

\begin{figure}
\begin{center}
\subfigure{\includegraphics[width = 3in]{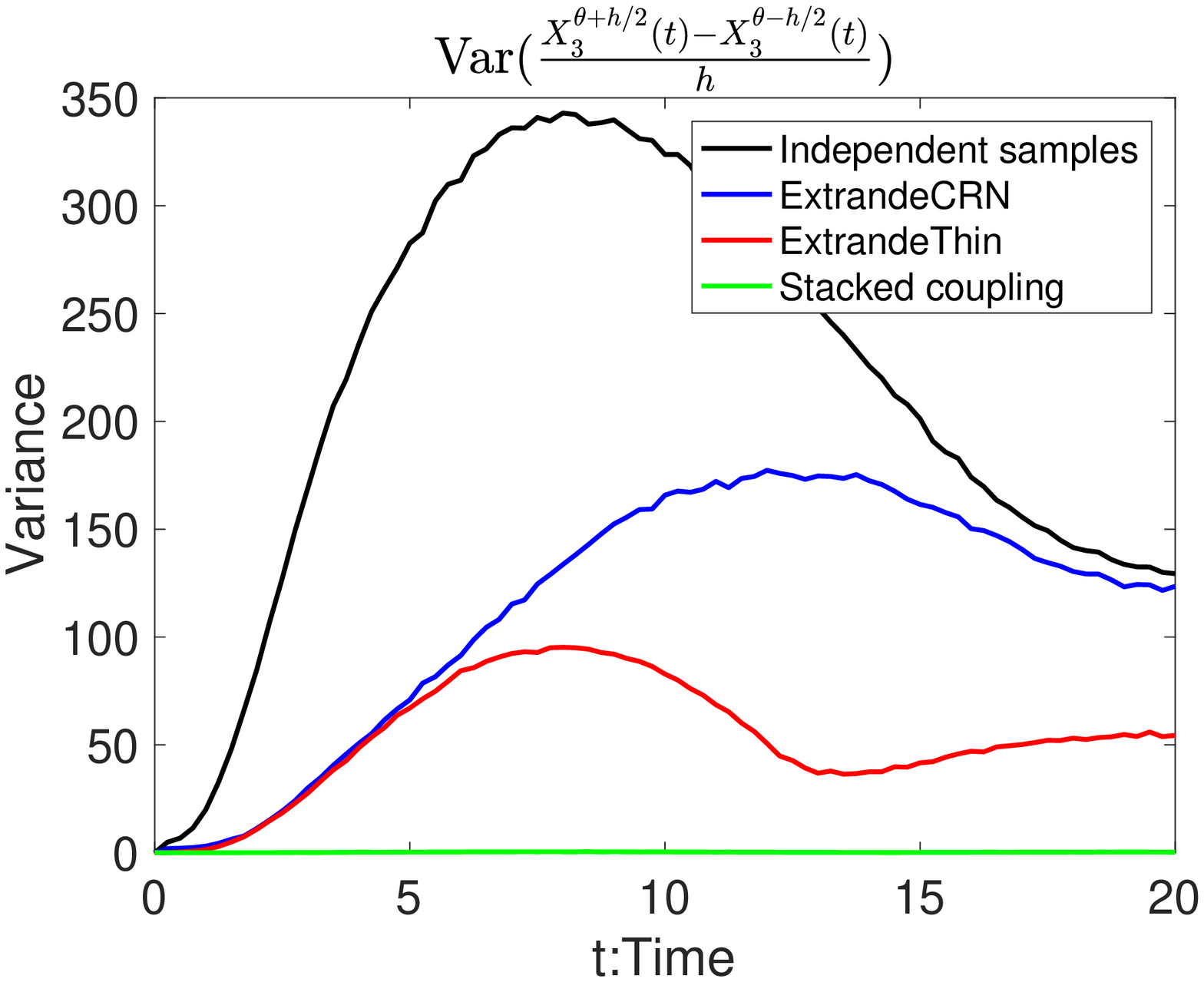}}
\subfigure{\includegraphics[width = 3in]{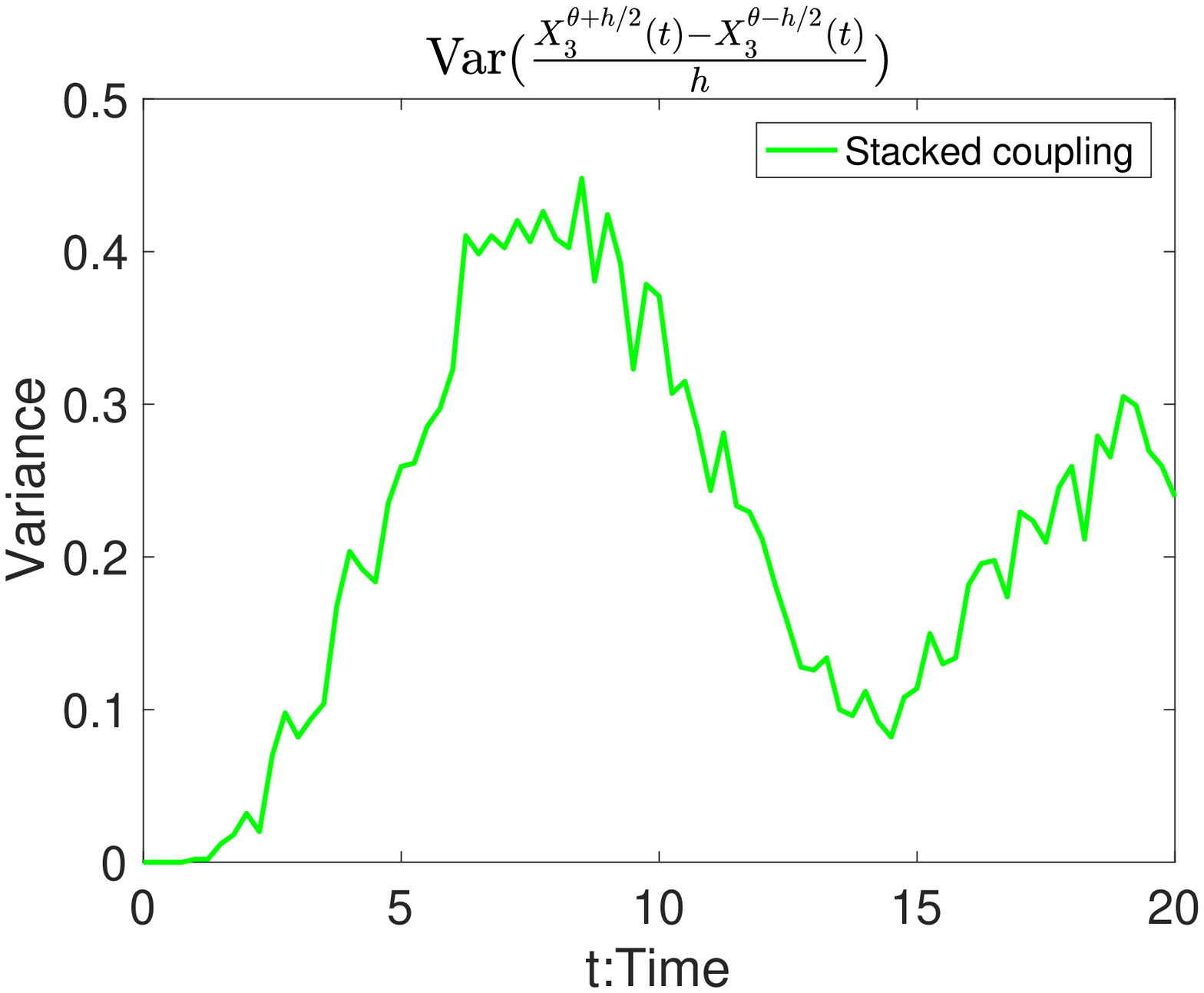}}
\end{center}
\caption{Estimates of $\text{Var}\big(\frac{X_{3}^{\theta+h/2}(t) -X_{3}^{\theta-h/2}(t))}{h}\big)$ at $\theta = 15$ for the different couplings  for the model in Example \ref{numexample1}.}
\label{fig:senmod4}
\end{figure}

\end{itemize}

Finally we will look at the long time behavior of these couplings. Analogously to \cite{AndCFD2012}, we expect that  as time gets large the  variance associated with ExtrandeCRN will converge to the variance associated with independent samples. However, we expect that the variance for the stacked coupling should remain small.   To test these hypotheses, we simulated the model until an end time of $T=200$ and plotted the resulting variances  in Figure \ref{fig:senmodL1}.  We see that this  simulation agrees with  our intuition. \hfill $\square$

\begin{figure}[!htbp]
\subfigure{\includegraphics[width = 3in]{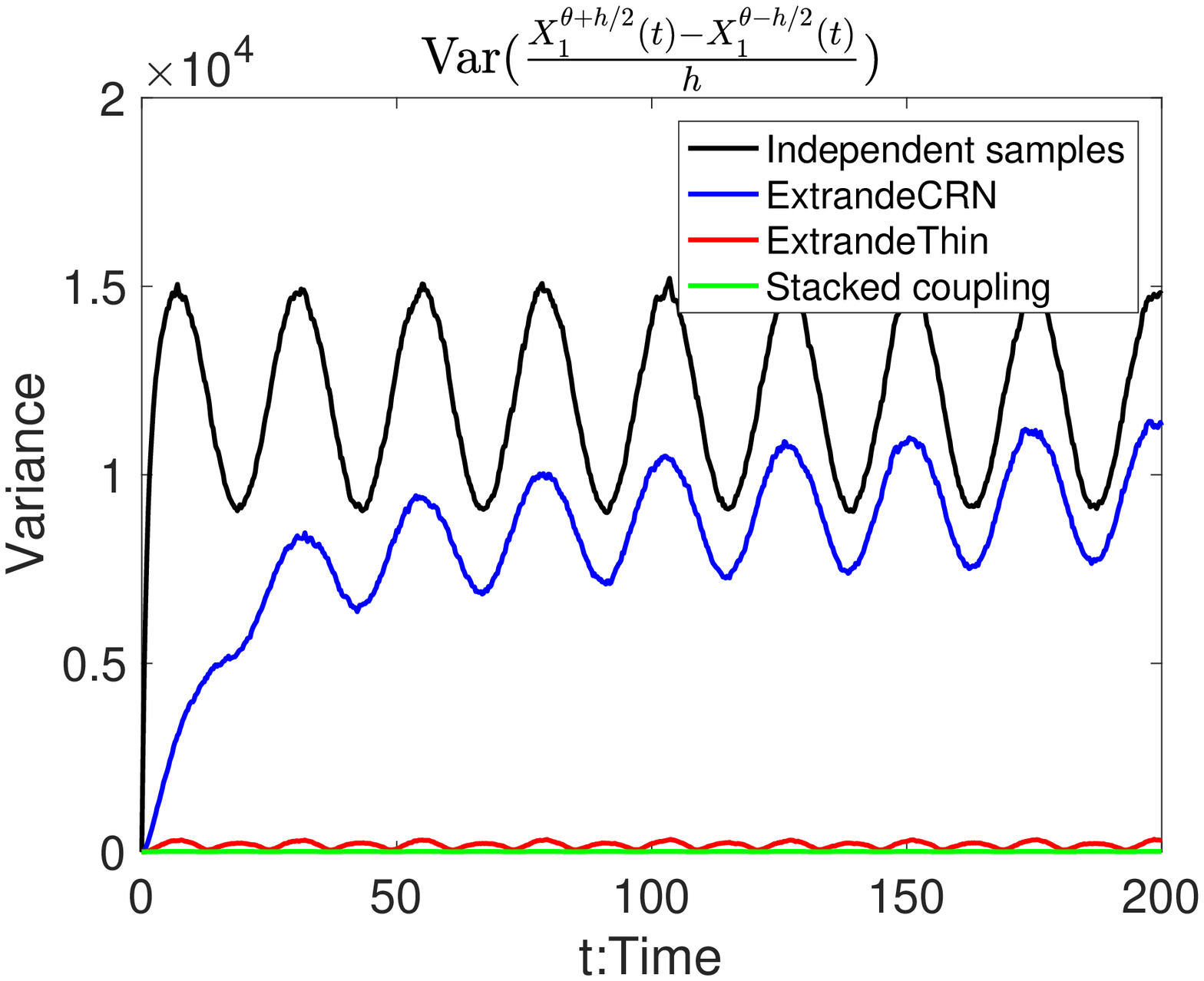}}
\subfigure{\includegraphics[width = 3in]{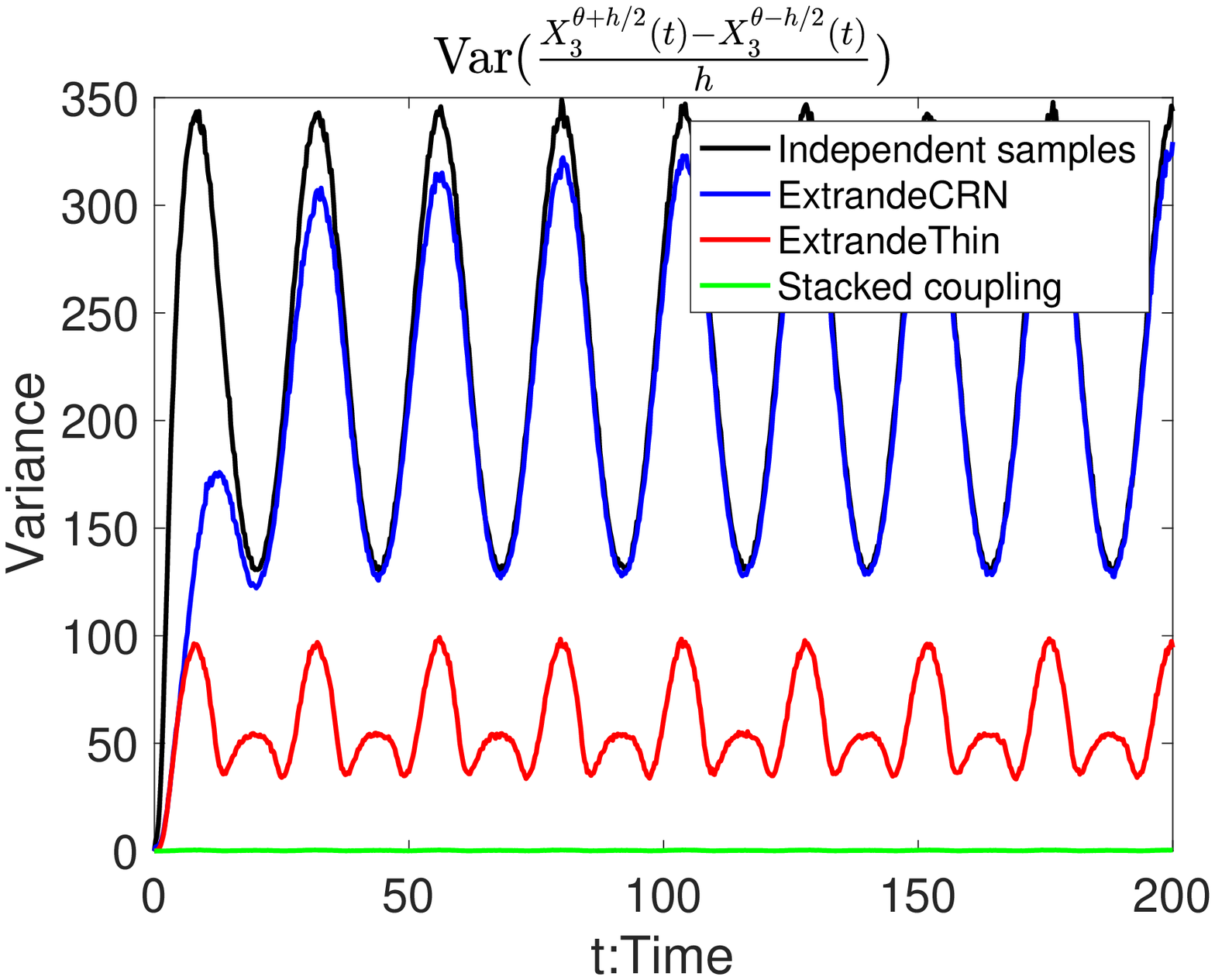}}
\caption{Estimates of $\text{Var}\big(\frac{X_{1}^{\theta+h/2}(t) -X_{1}^{\theta-h/2}(t)}{h}\big)$ and $\text{Var}\big(\frac{X_{3}^{\theta+h/2}(t) -X_{3}^{\theta-h/2}(t)}{h}\big)$ at $\theta = 15$ for the different couplings  for the model in Example \ref{numexample1} for $t \in [0,200]$.}
\label{fig:senmodL1}
\end{figure}

\end{example}

\begin{example}\label{numexample:2}
We continue to consider  Example  \ref{example:dimer}, except now we  study the effect of perturbing the degradation rate of mRNA in reaction $(R3)$ on the expected number of mRNA and dimers.  That is we parameterize the model by $\theta$ in $\lambda_3^\theta(x) = \theta x_1$ and consider
\[
	\frac{d}{d\theta}  E[X_{1}^\theta(t)]\big|_{\theta=1} \quad \text{and} \quad \frac{d}{d\theta}  E[X_{3}^\theta(t)]\big|_{\theta=1}.
\]
We use a perturbation of $h=0.05$ and, as before, 50,000 sample paths for each coupling.

We begin with the sensitivity for mRNA.
\begin{itemize}
\item In Figure \ref{fig:senmod5} we plot the estimated sensitivities of the different estimators as functions of time.
\begin{figure}
     \centering
      \includegraphics[width=3in]{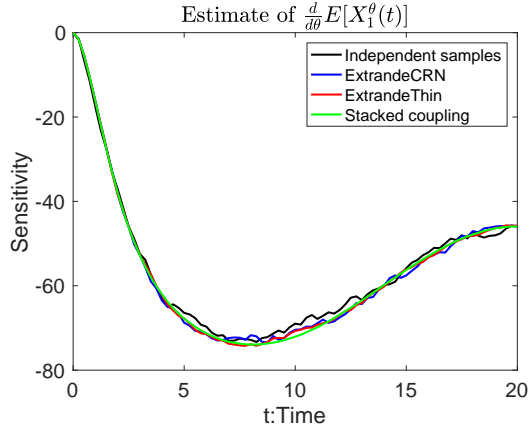} 
      \caption{Estimation of $\frac{d}{d\theta}E[X_{1}^\theta(t)]$ at $\theta = 1$ using finite differences and four different couplings for the model of Example \ref{numexample:2}.}
      \label{fig:senmod5}
\end{figure}

\item In Figure \ref{fig:senmod6} we plot the estimates of the variances
\[
	\text{Var}\bigg(\frac{X_{1}^{\theta+h/2}(t) -X_{1}^{\theta-h/2}(t))}{h}\bigg)
\]
at $\theta = 1$.
\begin{figure}
\subfigure{\includegraphics[width = 3in]{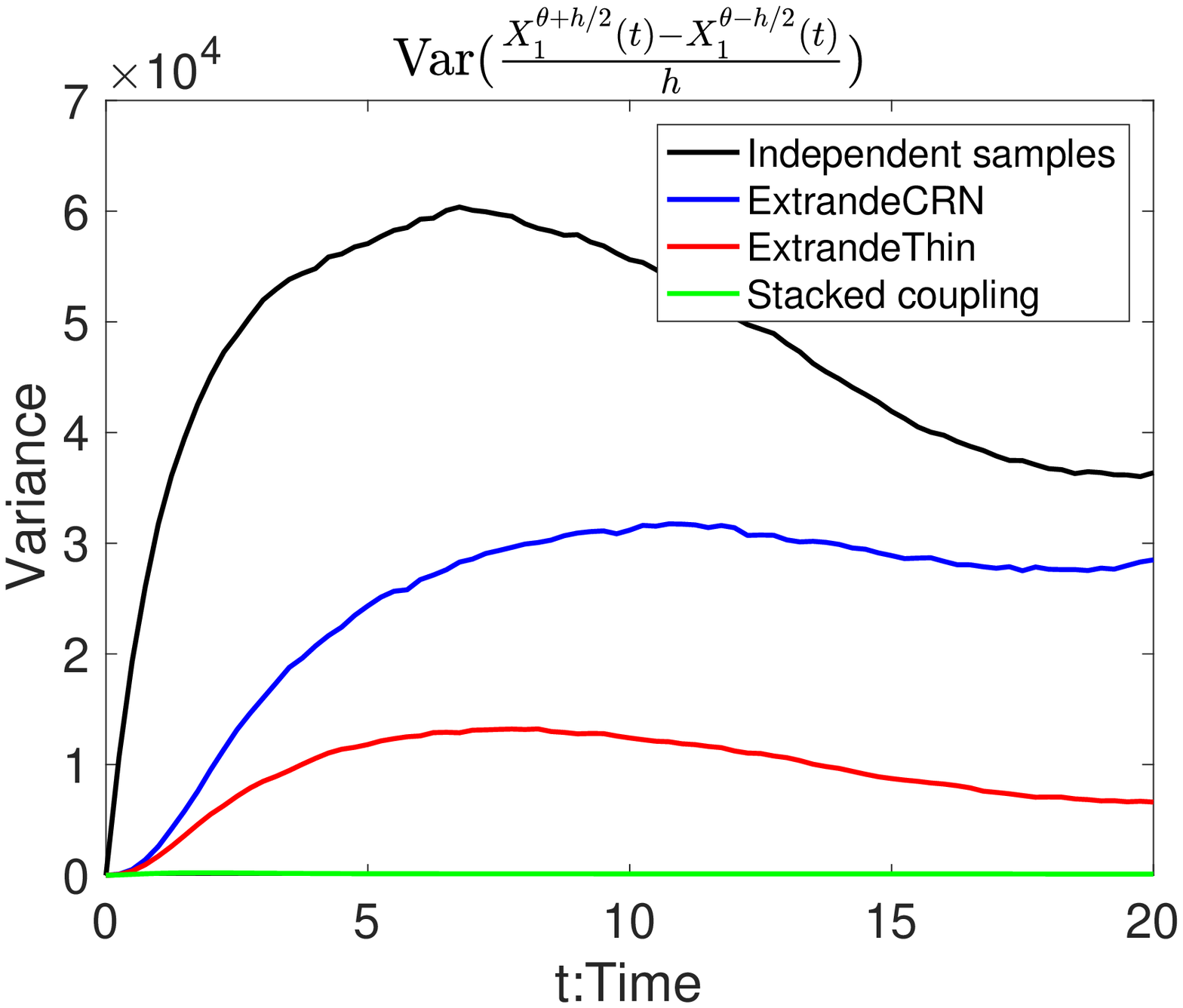}} \label{4(a)}
\subfigure{\includegraphics[width = 3in]{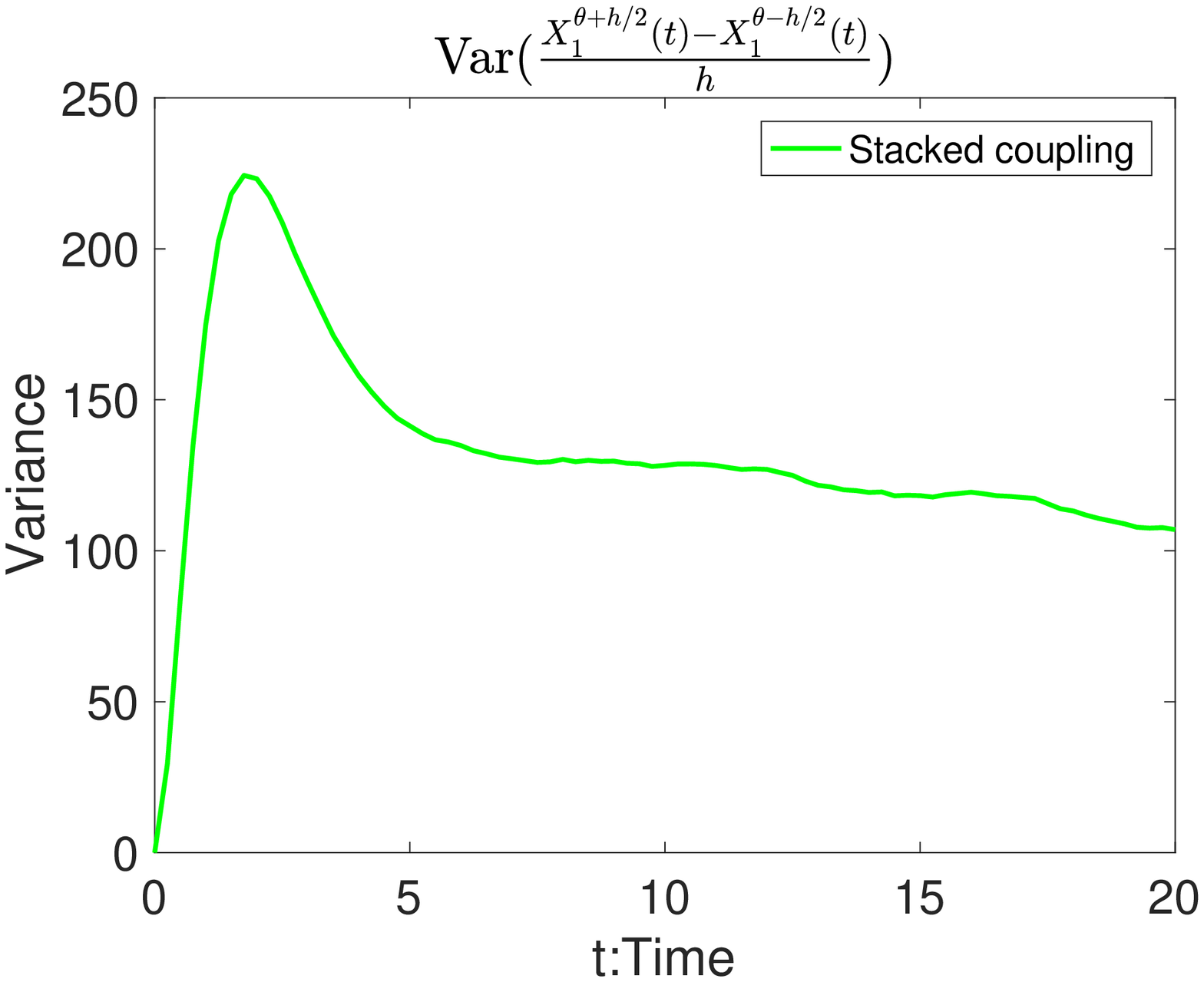}} \label{4(b)}
\caption{Estimates of $\text{Var}\big(\frac{X_{1}^{\theta+h/2}(t) -X_{1}^{\theta-h/2}(t))}{h}\big)$ at $\theta = 1$ for the different couplings  for the model in Example \ref{numexample:2}.}
\label{fig:senmod6}
\end{figure}
We once again observe that the stacked coupling produces an estimator with a variance that is dramatically lower than the other couplings.
\end{itemize}

We turn to the sensitivity for the dimers.

\begin{itemize}
\item In Figure \ref{fig:senmod7} we plot the estimated sensitivities of the different estimators as functions of time.
\begin{figure}
     \centering
      \includegraphics[width=3in]{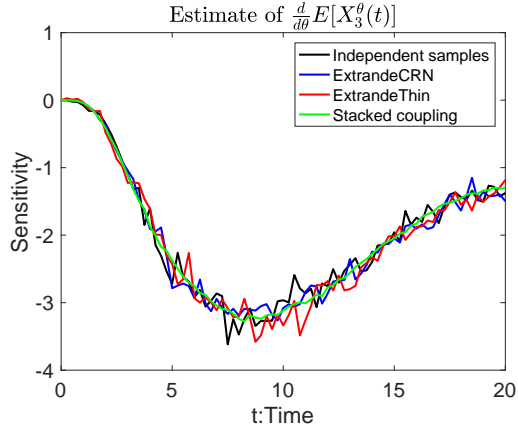} 
      \caption{Estimation of $\frac{d}{d\theta}E[X_{3}^\theta(t)]$ at $\theta = 1$ using finite differences and four different couplings for the model of Example \ref{numexample:2}.}
      \label{fig:senmod7}
\end{figure}

\item In Figure \ref{fig:senmod8} we plot the estimates of the variances
\[
	\text{Var}\bigg(\frac{X_{3}^{\theta+h/2}(t) -X_{3}^{\theta-h/2}(t))}{h}\bigg)
\]
at $\theta = 1$.
\begin{figure}[h!]
     \centering
      \includegraphics[width=3in]{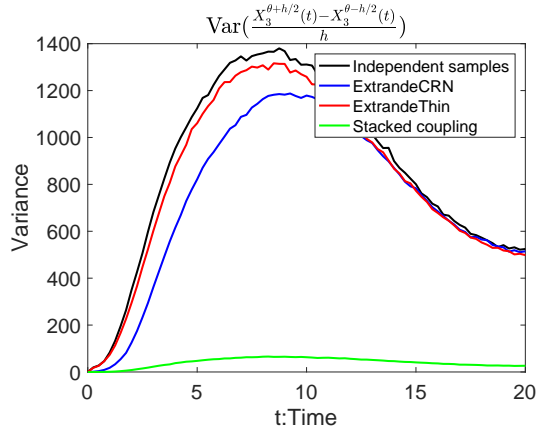} 
      \caption{Estimates of $\text{Var}\big(\frac{X_{3}^{\theta+h/2}(t) -X_{3}^{\theta-h/2}(t))}{h}\big)$ at $\theta = 1$ for the different couplings  for the model in Example \ref{numexample:2}.}
      \label{fig:senmod8}
\end{figure}
The stacked coupling once again provides a dramatically lower variance.
\end{itemize}

Finally,  we simulated the processes until time $T=200$ to see the long term behavior of these couplings.  The result is found in Figure \ref{fig:senmodL2}. Once again, the ExtrandeCRN coupling decoupled very quickly and the stacked coupling kept the variance small throughout the simulation. Interestingly,  in the right most image in Figure \ref{fig:senmodL2} we see that the Extrande thinning coupling did not perform well, as it decoupled almost immediately. \hfill $\square$

\begin{figure}[!htbp]
\subfigure{\includegraphics[width = 3in]{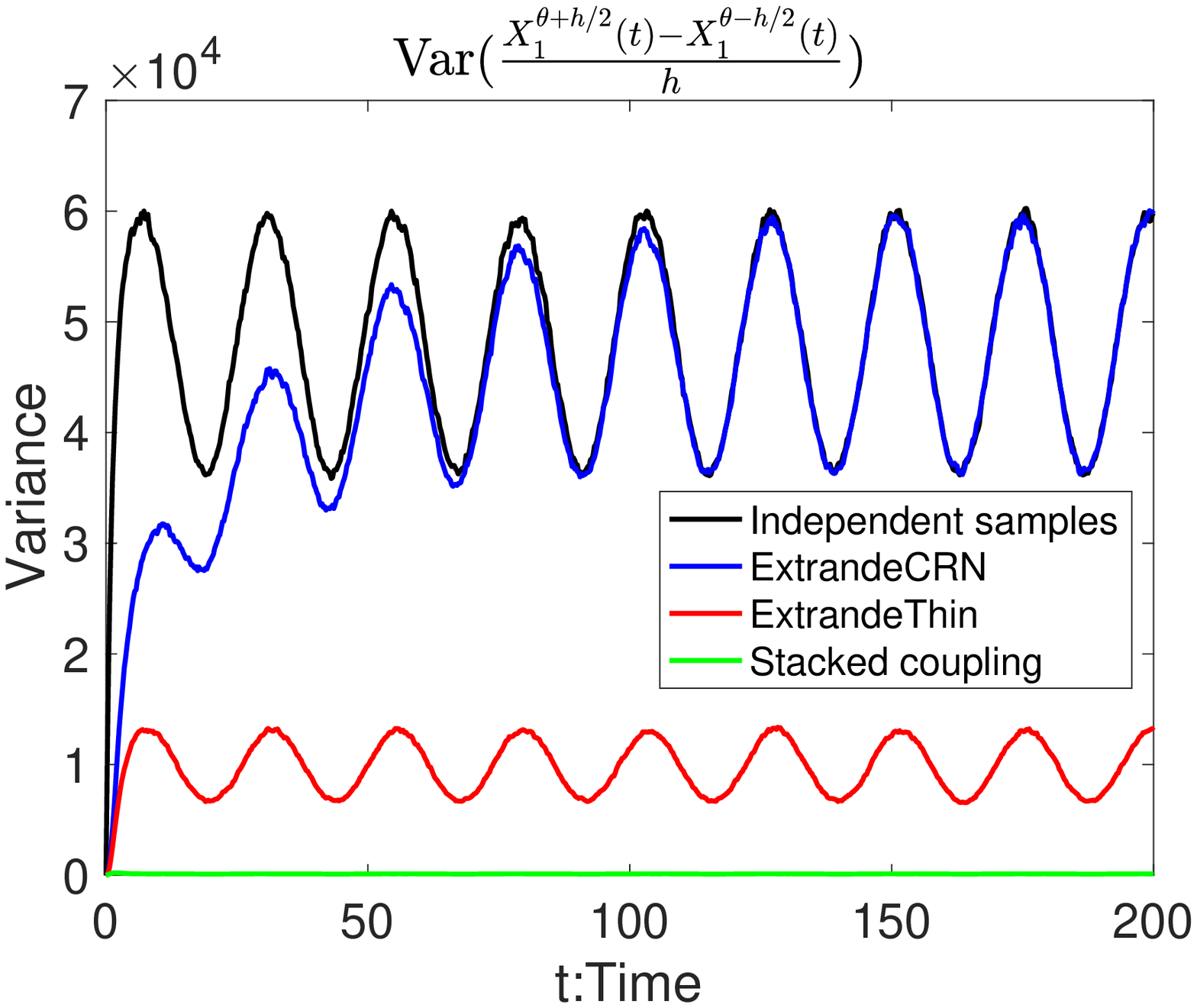}}
\subfigure{\includegraphics[width = 3in]{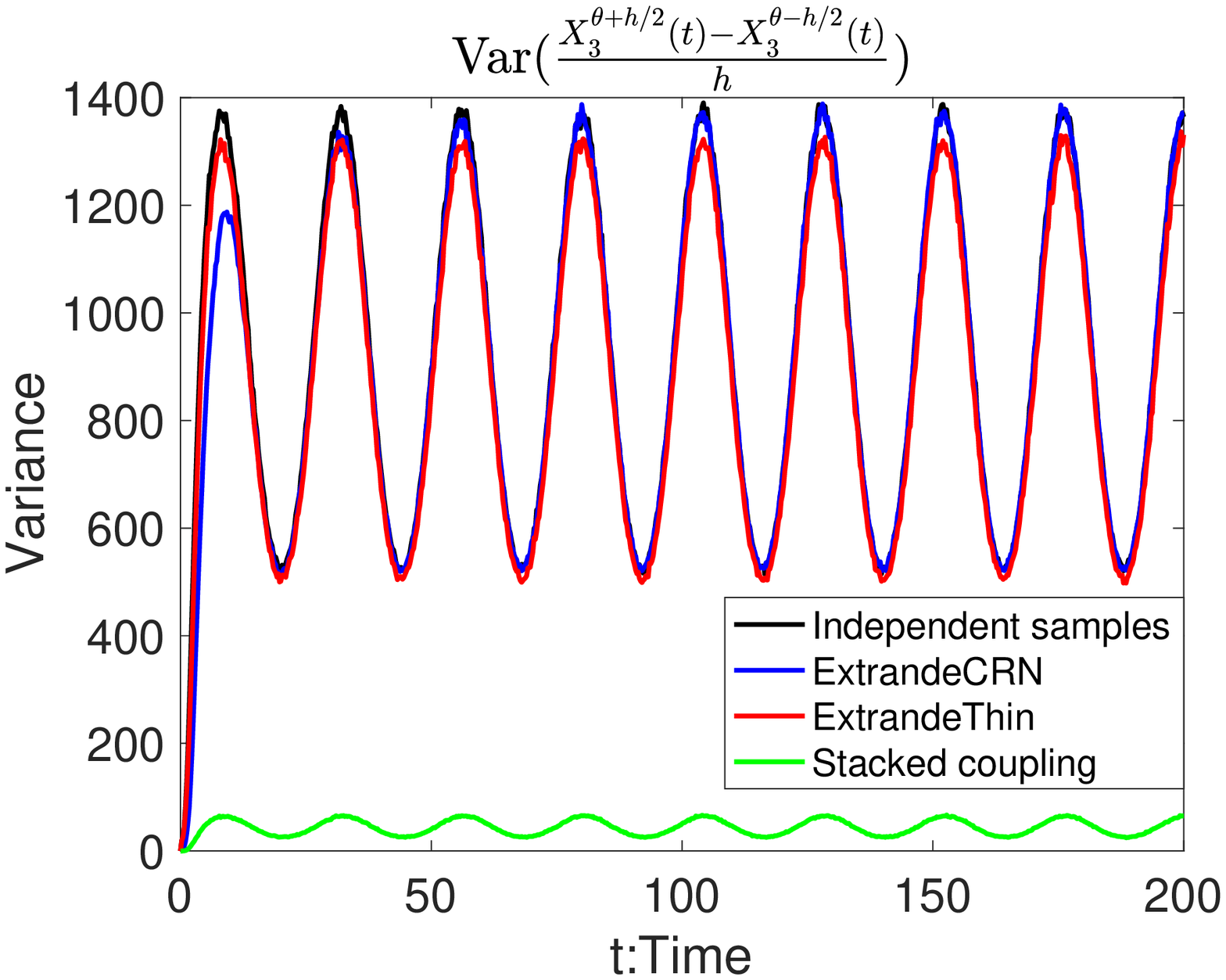}}
\caption{Estimates of $\text{Var}\big(\frac{X_{3}^{\theta+h/2}(t) -X_{3}^{\theta-h/2}(t)}{h}\big)$ and $\text{Var}\big(\frac{X_{3}^{\theta+h/2}(t) -X_{3}^{\theta-h/2}(t)}{h}\big)$ at $\theta = 15$ for the different couplings  for the model in Example \ref{numexample1}.}
\label{fig:senmodL2}
\end{figure}

\end{example}

\begin{example}
We turn to  Example Model \ref{example:SIR} of Section \ref{sec:couplings1}. As we recall, the quantity we are interested in is the sensitivity of the extinction probability, $p$, of pathogens within 10 years against all the parameters in our model. For each of the couplings, we utilized 5,000 sample paths and used a perturbation of $h=(\text{value of the parameter})\times0.05$.

Due to the number of results we are providing for this model, we provide our data in tabular form in Tables \ref{my-label1}, \ref{my-label2}, \ref{my-label3}, \ref{my-label4}, and \ref{my-label5}.  Each table returns the following information for each method:
\begin{enumerate}
\item An estimate of the sensitivity.
\item An estimate of $\text{Var}(\frac{p^{\theta+h/2}-p^{\theta-h/2}}{h})$.
\item The runtime in CPU seconds required to perform the calculation.

\item The number of random variables required by the method to perform the calculation.
\end{enumerate}
The key takeaway of the results is that the stacked coupling is always more efficient than the other couplings.  \hfill $\square$

\begin{table}[!htbp]
\centering
\begin{tabular}{|c|c|c|c|c|}
\hline
           & \textbf{Sensitivity} & \textbf{Variance}   & \textbf{runtime} & \# RVs \\ \hline
\textbf{Independent samples}  & 1.3640      &    199.0593     &  1932.5374   &  $1.3\times 10^{10}$            \\ \hline
\textbf{ExtrandeCRN}       & 0.8440             &  152.6782       & 1904.3838    & $1.3\times 10^{10}$              \\ \hline
\textbf{Extrande thinning} & 1.0400             &  193.8372        & 2386.6145   &  $9.2\times 10^{9}$               \\ \hline
\textbf{Stacked coupling}  & 0.8560             & 26.4726          & 1576.7281    &  $6.8\times 10^{10}$             \\ \hline
\end{tabular}
\caption{Sensitivity with respect to $m$ (death rate).}
\label{my-label1}
\end{table}

\begin{table}[!htbp]
\centering
\begin{tabular}{|c|c|c|c|c|}
\hline
           &  \textbf{Sensitivity}  & \textbf{Variance}  & \textbf{runtime} & \# RVs\\ \hline
\textbf{Independent samples}     &   0.2697   &    1.3348      & 1884.4148   &   $1.3\times 10^{10}$       \\ \hline
\textbf{ExtrandeCRN}                 &   0.2660   &    1.0150      & 1890.0392   &   $1.3\times 10^{10}$     \\ \hline
\textbf{Extrande thinning}           &   0.2577   &    0.7443      & 2015.2149   &   $7.6\times 10^{9}$      \\ \hline
\textbf{Stacked coupling}            &   0.2600   &    0.4603      & 1632.2105   &   $7.0\times 10^{9}$    \\ \hline
\end{tabular}
\caption{Sensitivity with respect to $\gamma$ (recovery rate).}
\label{my-label2}
\end{table}
\begin{table}[!htbp]
\centering
\begin{tabular}{|c|c|c|c|c|}
\hline
           &  \textbf{Sensitivity}  & \textbf{Variance} &  \textbf{runtime} & \# RVs\\ \hline
\textbf{Independent samples}      & -0.5882     &    742.4462    & 1901.3166     &  $1.3\times 10^{10}$ \\ \hline
\textbf{ExtrandeCRN}                  & -0.1146     &    592.1777    & 1899.8093     & $1.3\times 10^{10}$  \\ \hline
\textbf{Extrande thinning}            & 0.2368      &    589.2161    & 2166.4572     & $8.3\times 10^{9}$   \\ \hline
\textbf{Stacked coupling}             & -0.5118     &   136.6219     & 1559.3501     & $6.7\times 10^{9}$  \\ \hline
\end{tabular}
\caption{Sensitivity with respect to $\phi$ (phase parameter).}
\label{my-label3}
\end{table}
\begin{table}[!htbp]
\centering
\begin{tabular}{|c|c|c|c|c|}
\hline
           &  \textbf{Sensitivity}  & \textbf{Variance} & \textbf{runtime} &\# RVs \\ \hline
\textbf{Independent samples}     & 0.3140    &   12.6939    & 1876.2869    &   $1.3\times 10^{10}$  \\ \hline
\textbf{ExtrandeCRN}       		 & 0.4630    &   9.1225      & 1897.8293    &   $1.3\times 10^{10}$  \\ \hline
\textbf{Extrande thinning}           & 0.3790    &   5.9826      & 2015.0119     &   $7.5\times 10^{9}$     \\ \hline
\textbf{Stacked coupling}            & 0.4010    &   2.9548      & 1601.5408    &   $6.9\times 10^{9}$     \\ \hline
\end{tabular}
\caption{Sensitivity with respect to $R_0$ (reproductive ratio).}
\label{my-label4}
\end{table}
\begin{table}[!htbp]
\centering
\begin{tabular}{|c|c|c|c|c|}
\hline
           &  \textbf{Sensitivity} & \textbf{Variance} &  \textbf{runtime} & \# RVs  \\ \hline
\textbf{Independent samples}     & 0.0276    &   1.9540    & 1947.5821   &  $1.3\times 10^{10}$  \\ \hline
\textbf{ExtrandeCRN}      		 & 0.0292    &    1.9715   & 2024.7114   &  $1.3\times 10^{10}$ \\ \hline
\textbf{Extrande thinning} 		 & 0.0556    &    1.7997   & 2399.2836   &  $8.6\times 10^{10}$ \\ \hline
\textbf{Stacked coupling}  		 & 0.0140    &     0.1222  & 1645.1277   &  $6.5\times 10^{10}$ \\ \hline
\end{tabular}
\caption{Sensitivity with respect to $s$ (synchrony parameter).}
\label{my-label5}
\end{table}
\end{example}

\pagebreak

\subsection{Numerical examples for multilevel Monte Carlo}
\label{subsec: MLMC}

\begin{example} We  consider the Markov modulated process as described in Example Model \ref{example:MMP}.  We will estimate the expected  number of molecules at  time $T=2$, specifically we estimate
\[
	E[X_{i}(2)],\quad \text{ for each of  $i=1,2,3,4$}.
\]
The upper bound  $\overline\lambda_1$ is taken to be $5/1000$ as $k_1(t)$ is always less than or equal to 5.  For each of the MLMC methods, we used $M=4$ with a  time step of $4^{-2}$ at the coarsest level and a time step of $4^{-3}$ at the finest level.  

We compare direct simulation using the Extrande method (Algorithm \ref{algorithm1}) and multilevel Monte Carlo with the Stacked coupling at the finest level.  We simulated each method until the estimator standard deviation fell below a tolerance of $\varepsilon = 0.1$.  In Table   \ref{table:mlmccor1}, the estimates for the expectations are provided.  In Table \ref{table:mlmccot1}, we report the CPU time required by the different methods for the different expectations.  Finally, in Table \ref{table:mlmccorv1} we report the number of random variables utilized by the different methods for the different expectations.  \hfill $\square$

\begin{table}[htbp!]
\centering\begin{tabular}{|c || c | c | c | c | c || }
 \hline
 \multicolumn{5}{|c|}{Comparison of estimated expectation} \\
 \hline
 Method & $E[X_{1}(2)]$ &$E[X_{2}(2)]$ &$E[X_{3}(2)]$&$E[X_{4}(2)]$ \\
 \hline
Extrande & 335.5 & 768.2 & 231.9 & 432.5 \\
\hline
MLMC & 335.6 & 768.0 &  231.9 & 432.5 \\
 \hline
\end{tabular}
\caption{Estimates for desired expectations.}
\label{table:mlmccor1}
\end{table}

\begin{table}[htbp!]
\centering\begin{tabular}{|c || c | c | c | c | c || }
 \hline
 \multicolumn{5}{|c|}{Comparison of time cost (CPU seconds)} \\
 \hline
 Method & $E[X_{1}(2)]$ &$E[X_{2}(2)]$ &$E[X_{3}(2)]$&$E[X_{4}(2)]$ \\
 \hline
Extrande & $2.91\times 10^{3}$  & $6.82\times 10^{2}$    & $6.66\times 10^{2}$ & $1.98\times 10^{3}$ \\
\hline
MLMC & $4.59\times 10^{2}$ & $1.35\times 10^{2}$ & $1.31\times 10^{2}$ & $3.04\times 10^{2}$ \\
 \hline
\end{tabular}
\caption{CPU times required.}
\label{table:mlmccot1}
\end{table}
 
\begin{table}[htbp!]
\centering\begin{tabular}{|c || c | c | c | c | c || }
 \hline
 \multicolumn{5}{|c|}{Comparison of random variable generated} \\   \hline
 Method & $E[X_{1}(2)]$ &$E[X_{2}(2)]$ &$E[X_{3}(2)]$&$E[X_{4}(2)]$ \\ \hline
Extrande & $3.1\times 10^{10}$  & $7.2\times 10^{9}$ & $7.2\times 10^{9}$  & $2.0\times 10^{10}$ \\ \hline
MLMC & $6\times 10^{8}$ & $2.1\times 10^{8}$ &  $2.1\times 10^{8}$ & $3.3\times 10^{8}$  \\ \hline
\end{tabular}
\caption{Number of random variables required.}
\label{table:mlmccorv1}
\end{table}
\end{example}

\begin{example}  
\label{example:mlmc2}
We return to  the mRNA-protein-dimer model of Example  \ref{example:dimer}.  Depending upon the expectation being calculated,  a standard implementation of MLMC will not perform well for this model.  The reason for this is that while some of the propensity functions are high (such as reaction (R4)), others are quite low.  In particular, the rates of reactions (R5) and (R6) will remain at approximately 1 throughout the computation.   We therefore choose a different approximate model than the usual Euler tau-leap model.  In particular, we take our approximate process $Z_{\ell}$ to have a time discretization parameter of $h_\ell$ and
\begin{itemize}
\item have constant propensity functions for reactions (R1) -- (R5) over time steps of size $h_\ell$ (similar to standard tau-leaping), and
\item have intensity function $\lambda_6(Z_{\ell}(s)) = 10 Z_{\ell,3}(s)$ for all $s \ge 0$.
\end{itemize}
Thus, we are not using an Euler discretization for the sixth reaction channel, but are doing so for the first five reaction channels.  We then couple the relevant processes at each level via the stacked coupling.

We compare three different methods for the estimation of each of $E[X_i(20)]$ for $i \in \{1,2,3\}$: standard Monte Carlo with the Extrande method, standard MLMC (i.e.~with propensity function $\lambda_6$ also being discretized), and MLMC using no discretization for $\lambda_6$.  We will call the last method MLMC6.   We simulated each method until the estimator standard deviation fell below a tolerance of 0.05, 5, and 0.005 for $E[X_1(20)]$, $E[X_2(20)]$, and $E[X_3(20)]$, respectively. 
For each of the MLMC methods, we used $M=4$ with a  time step of $4^{-1}$ at the coarsest level and a time step of $4^{-2}$ at the finest level.    

 In Table   \ref{table:mlmccor2}, the estimates for the expectations are provided.  In Table \ref{table:mlmccot2} we report the CPU times required by the different methods for the different expectations.  Finally, in Table \ref{table:mlmccorv2} we report the number of random variables utilized by the different methods for the different expectations.

\begin{table}[htbp!]
\centering\begin{tabular}{|c || c | c | c | c | c || }
 \hline
 \multicolumn{4}{|c|}{Comparison of results} \\
 \hline
 Method & $E[X_{1}(20)]$ &$E[X_{2}(20)]$ &$E[X_{3}(20)]$ \\
 \hline
Extrande & 46.01 & $4.58\times 10^{3}$ & 0.6413 \\
\hline
MLMC & 45.94 & $4.59\times 10^{3}$ &  0.6329\\
 \hline
MLMC6 & 46.01 & $4.58\times 10^{3}$ &  0.6393\\
 \hline
\end{tabular}
\caption{Estimates for desired expectations.}
\label{table:mlmccor2}
\end{table}

\begin{table}[htbp!]
\centering\begin{tabular}{|c || c | c | c | c | c || }
 \hline
 \multicolumn{4}{|c|}{Comparison of time cost(second)} \\
 \hline
  Method & $E[X_{1}(20)]$ &$E[X_{2}(20)]$ &$E[X_{3}(20)]$ \\
 \hline
Extrande & $4.56\times 10^{3}$ & $2.22\times 10^{3}$ & $7.40\times 10^{3}$ \\
\hline
MLMC & $1.09\times 10^{2}$ & 45.18 & $2.36\times 10^{3}$ \\
 \hline
  MLMC6 & $1.11\times 10^{2}$ & 57.94  & $2.33\times 10^{2}$ \\
 \hline
\end{tabular}
\caption{CPU times required.}
\label{table:mlmccot2}
\end{table}

\begin{table}[htbp!]
\centering\begin{tabular}{|c || c | c | c | c | c || }
 \hline
 \multicolumn{4}{|c|}{Comparison of random variable generated} \\
 \hline
 Method & $E[X_{1}(20)]$ &$E[X_{2}(20)]$ &$E[X_{3}(20)]$ \\
 \hline
Extrande & $8.7\times 10^{9}$ & $4.7\times 10^{9}$ & $1.7\times 10^{10}$ \\
\hline
MLMC & $3.5\times 10^{8}$ & $1.4\times 10^{8}$ &  $8.2\times 10^{9}$ \\
 \hline
MLMC6 & $3.2\times 10^{8}$ & $1.6\times 10^{8}$ & $7.2\times 10^{8}$ \\
 \hline
\end{tabular}
\caption{Number of random variables required.}
\label{table:mlmccorv2}
\end{table}
We can see from the tables that both MLMC and MLMC6 provide improvements over direct Monte Carlo with the Extrande method and have similar performances when calculating the expected number of mRNA and protein molecules.  However, they differed by a factor of 10 in performance when calculating the expected number of dimers.  This  demonstrates that clever choices of approximate processes can have  a dramatic effect on your estimator performance.\hfill $\square$
\end{example}

\section{Conclusion}
\label{sec:discussion}
There have been a number of papers recently on the efficient simulation of biochemical processes with time dependent parameters.  In this paper, we studied how to efficiently couple such processes in order to efficiently perform parametric sensitivity analysis and multilevel Monte Carlo.  
In particular, through different mathematical representations for the models of interest, we developed three non-trivial coupling strategies.  Through a number of examples, we then demonstrated that the \textit{stacked coupling} appears to be the most efficient. 
  Future work will focus on providing analytical results related to the variances of the different coupling strategies.

 \bibliographystyle{plain}
\bibliography{Thinning}

\begin{thebibliography}{10}

\bibitem{anderson2007modified}
David~F. Anderson.
\newblock A modified next reaction method for simulating chemical systems with
  time dependent propensities and delays.
\newblock {\em The Journal of chemical physics}, 127(21):214107, 2007.

\bibitem{Anderson2007b}
David~F. Anderson.
\newblock Incorporating postleap checks in tau-leaping.
\newblock {\em J. Chem. Phys.}, 128(5):054103, 2008.

\bibitem{AndCFD2012}
David~F. Anderson.
\newblock An efficient finite difference method for parameter sensitivities of
  continuous time {M}arkov chains.
\newblock {\em SIAM Journal on Numerical Analysis}, 50(5):2237 -- 2258, 2012.

\bibitem{ACK2016}
David~F. Anderson, Daniele Cappelletti, and Thomas~G. Kurtz.
\newblock Finite time distributions of stochastically modeled chemical systems
  with absolute concentration robustness.
\newblock {\em SIAM J. Appl. Dyn. Syst.}, 16(3):1309--1339, 2017.

\bibitem{AEJ2014}
David~F. Anderson, Germ\'an~A. Enciso, and Matthew~D. Johnston.
\newblock Stochastic analysis of biochemical reaction networks with absolute
  concentration robustness.
\newblock {\em Royal Society Interface}, 11:20130943, 2014.

\bibitem{AET2015}
David~F. Anderson, Bard Ermentrout, , and Peter~J. Thomas.
\newblock Stochastic representations of ion channel kinetics and exact
  stochastic simulation of neuronal dynamics.
\newblock {\em J. Comp. Neuro.}, 38(1):67--82, 2015.

\bibitem{AndersonGangulyKurtz}
David~F. Anderson, Arnab Ganguly, and Thomas~G. Kurtz.
\newblock Error analysis of tau-leap simulation methods.
\newblock {\em Annals of Applied Probability}, 21(6):2226 -- 2262, 2011.

\bibitem{AndHigham2012}
David~F. Anderson and Desmond~J. Higham.
\newblock Multi-level {M}onte {C}arlo for continuous time {M}arkov chains, with
  applications in biochemical kinetics.
\newblock {\em SIAM: Multiscale Modeling and Simulation}, 10(1):146 -- 179,
  2012.

\bibitem{AHS2014}
David~F. Anderson, Desmond~J. Higham, and Yu~Sun.
\newblock Complexity of multilevel {M}onte {C}arlo tau-leaping.
\newblock {\em SIAM J. Numer. Anal.}, 52(6):3106--3127, 2014.

\bibitem{AHS2016}
David~F. Anderson, Desmond~J. Higham, and Yu~Sun.
\newblock Multilevel {M}onte {C}arlo for stochastic differential equations with
  small noise.
\newblock {\em SIAM J. Numer. Anal.}, 54(2):505--529, 2016.

\bibitem{AHS2017}
David~F. Anderson, Desmond~J. Higham, and Yu~Sun.
\newblock Computational complexity analysis for {M}onte {C}arlo approximations
  of classically scaled population processes.
\newblock Submitted, available on arXiv:
  \url{https://arxiv.org/abs/1512.01588}, 2017.

\bibitem{AK2011}
David~F. Anderson and Thomas~G. Kurtz.
\newblock {\em Continuous time Markov chain models for chemical reaction
  networks}, chapter 1 in Design and Analysis of Biomolecular Circuits:
  Engineering Approaches to Systems and Synthetic Biology.
\newblock Springer, 2011.

\bibitem{AK2015}
David~F. Anderson and Thomas~G. Kurtz.
\newblock {\em Stochastic analysis of biochemical systems}.
\newblock Springer, 2015.

\bibitem{AndSkubak}
David~F. Anderson and Elizabeth~Skubak Wolf.
\newblock A finite difference method for estimating second order parameter
  sensitivities of discrete stochastic chemical reaction networks.
\newblock {\em J. Chem. Phys.}, 137(22):224112, 2012.

\bibitem{GlynnAsmussen2007}
S{\o}ren Asmussen and Peter~W. Glynn.
\newblock {\em Stochastic Simulation: Algorithms and Analysis}.
\newblock Stochastic modelling and applied probability. Springer, New York,
  2007.

\bibitem{duncan2015hybrid}
Andrew Duncan, Radek Erban, and Konstantinos Zygalakis.
\newblock Hybrid framework for the simulation of stochastic chemical kinetics.
\newblock {\em Journal of Computational Physics}, 326(1):398--419, 2016.

\bibitem{Gibson2000}
M.A. Gibson and J.~Bruck.
\newblock Efficient exact stochastic simulation of chemical systems with many
  species and many channels.
\newblock {\em J. Phys. Chem. A}, 105:1876--1889, 2000.

\bibitem{Giles2008}
Mike~B. Giles.
\newblock Multilevel {M}onte {C}arlo path simulation.
\newblock {\em Operations Research}, 56:607--617, 2008.

\bibitem{Gill2001}
D.~T. Gillespie.
\newblock Approximate accelerated simulation of chemically reaction systems.
\newblock {\em J. Chem. Phys.}, 115(4):1716--1733, 2001.

\bibitem{gillespie1976general}
Daniel~T. Gillespie.
\newblock A general method for numerically simulating the stochastic time
  evolution of coupled chemical reactions.
\newblock {\em Journal of computational physics}, 22(4):403--434, 1976.

\bibitem{gillespie1977exact}
Daniel~T. Gillespie.
\newblock Exact stochastic simulation of coupled chemical reactions.
\newblock {\em The journal of physical chemistry}, 81(25):2340--2361, 1977.

\bibitem{GuptaKhammash}
Ankit Gupta and Mustafa Khammash.
\newblock Unbiased estimation of parameter sensitivities for stochastic
  chemical reaction networks.
\newblock {\em SIAM Journal on Scientific Computing}, 35(6):A2598--A2620, 2013.

\bibitem{GK2014b}
Ankit Gupta and Mustafa Khammash.
\newblock An efficient and unbiased method for sensitivity analysis of
  stochastic reaction networks.
\newblock {\em Royal Society Interface}, 11(101):20140979, 2014.

\bibitem{GK2014}
Ankit Gupta and Mustafa Khammash.
\newblock Sensitivity analysis for stochastic chemical reaction networks with
  multiple time-scales.
\newblock {\em Electronic Journal of Probability}, 19(59):1--53, 2014.

\bibitem{Kurtz80}
Thomas~G. Kurtz.
\newblock Representations of markov processes as multiparameter time changes.
\newblock {\em Ann. Prob.}, 8(4):682--715, 1980.

\bibitem{lemaire2016exact}
Vincent Lemaire, Mich{\`e}le Thieullen, and Nicolas Thomas.
\newblock Exact simulation of the jump times of a class of piecewise
  deterministic markov processes.
\newblock arXiv preprint arXiv:1602.07871, 2016.

\bibitem{lester2016extending}
Christopher Lester, Ruth~E Baker, Michael~B Giles, and Christian~A Yates.
\newblock Extending the multi-level method for the simulation of stochastic
  biological systems.
\newblock {\em Bulletin of mathematical biology}, 78(8):1640--1677, 2016.

\bibitem{li2007analysis}
Tiejun Li.
\newblock Analysis of explicit tau-leaping schemes for simulating chemically
  reacting systems.
\newblock {\em Multiscale Modeling \& Simulation}, 6(2):417--436, 2007.

\bibitem{marchetti2016hrssa}
Luca Marchetti, Corrado Priami, and Vo~Hong Thanh.
\newblock Hrssa--efficient hybrid stochastic simulation for spatially
  homogeneous biochemical reaction networks.
\newblock {\em Journal of Computational Physics}, 317:301--317, 2016.

\bibitem{MTV2016a}
A.~Moraes, R.~Tempone, and P.~Vilanova.
\newblock Multilevel hybrid {C}hernoff tau-leap.
\newblock {\em BIT Numerical Mathematics}, 56(1):189--239, 2016.

\bibitem{peel2014effect}
Alison~J. Peel, J.R.C. Pulliam, A.D. Luis, R.K. Plowright, T.J. O'Shea, D.T.S.
  Hayman, J.L.N. Wood, C.T. Webb, and O.~Restif.
\newblock The effect of seasonal birth pulses on pathogen persistence in wild
  mammal populations.
\newblock {\em Proceedings of the Royal Society of London B: Biological
  Sciences}, 281(1786):20132962, 2014.

\bibitem{Khammash2010}
Muruhan Rathinam, Patrick~W. Sheppard, and Mustafa Khammash.
\newblock Efficient computation of parameter sensitivities of discrete
  stochastic chemical reaction networks.
\newblock {\em Journal of Chemical Physics}, 132:034103, 2010.

\bibitem{SAR2013}
Rishi Srivastava, David~F. Anderson, and James~B. Rawlings.
\newblock Comparison of finite difference based methods to obtain sensitivities
  of stochastic chemical kinetic models.
\newblock {\em Journal of Chemical Physics}, 138:074110, 2013.

\bibitem{thanh2015simulation}
Vo~Hong Thanh and Corrado Priami.
\newblock Simulation of biochemical reactions with time-dependent rates by the
  rejection-based algorithm.
\newblock {\em The Journal of chemical physics}, 143(5):054104, 2015.

\bibitem{voliotis2016stochastic}
Margaritis Voliotis, Philipp Thomas, Ramon Grima, and Clive~G Bowsher.
\newblock Stochastic simulation of biomolecular networks in dynamic
  environments.
\newblock {\em PLoS Comput Biol}, 12(6):e1004923, 2016.

\end{thebibliography}

\end{document}